\documentclass[11pt]{article}
\usepackage{amsmath,amssymb,eucal}
\usepackage[dvips]{graphicx}

\newtheorem{thm}{\bf Theorem}[section]
\newtheorem{cor}[thm]{\bf Corollary}
\newtheorem{lem}[thm]{\bf Lemma}

\newtheorem{prop}[thm]{\bf Proposition}
\newtheorem{ex}[thm]{\bf Example}
\newtheorem{defn}[thm]{\bf Definition}

\newtheorem{quest}[thm]{\bf Question}

\def\qed{$\Box$}

\begin{document}

\title{
Lens space surgeries on A'Campo's divide knots
}
\author{Yuichi YAMADA}
\date{}
\footnotetext[0]{%
}
\footnotetext[1]{%
This work was partially supported by Grant-in-Aid for
Scientific Research No.18740029, 
Japan Society for the Promotion of Science. 
} 
\maketitle
\begin{center}
Dedicated to Professor Takao Matumoto on the occasion of 
his 60th birthday.\end{center}
\begin{abstract}{
It is proved that 
every knot in the major subfamilies of J.~Berge's lens space surgery
(i.e., knots yielding a lens space by Dehn surgery)
is presented by an L-shaped (real) plane curve as 
a {\it divide knot} defined by N.~A'Campo in the context of 
singularity theory of complex curves.
For each knot given by Berge's parameters,
the corresponding plane curve is constructed.
The surgery coefficients are also considered.
Such presentations support us to study 
each knot itself, 
and the relationship among the knots in the set of lens space surgeries.
}\end{abstract}

\section{Introduction}~\label{sec:intro}
If $r/s$ Dehn surgery on a knot $K$ in $S^3$ yields the
lens space $L(p,q)$, we call the pair $(K, r/s)$ a {\it lens space surgery},
and we also say that $K$ admits a lens space surgery,
and that $r/s$ is the {\it coefficient} of the lens space surgery.
The task of classifying lens space surgeries, especially 
knots that admit lens space surgeries has been
a focal point in low-dimensional topology
and has been invigorated of late by results from the
Heegaard Floer homology theories of Ozsv\'{a}th--Sz\'{a}bo
\cite{OSz} (see also \cite{He}, \cite{Ta} and so on).
Before the first hyperbolic examples
found by Fintushel--Stern \cite{FS} in 1980, only
torus knots (Moser \cite{Mo})
and their $2$-cables (Bailey--Rolfsen \cite{BR}) were known.
After \cite{FS}, some more examples were found (see \cite{Ma}).
In 1990, Berge \cite{Bg}
pointed out a \lq\lq mechanism\rq\rq \ %
of known lens space surgery,
that is, {\it doubly-primitive knots} in the Heegaard surface of genus $2$.
Berge also gave a conjecturally complete list of such knots,
described them by Osborne--Stevens's
\lq\lq R-R diagrams\rq\rq \ %
in \cite{OSt}, and classified such knots into three families,
and into 12 types in detail:
\begin{enumerate}
\item[(1)] {\it Knots in a solid torus} (Type I, II, ... and VI)
\par
Dehn surgery along a knot in a solid torus whose resulting
manifold is also a solid torus.
This family was studied in \cite{Bg2}.
\item[(2)] {\it Knots in genus-one fiber surface} (Type VII and VIII)
\par
Dehn surgery along a knot in the genus-one fiber surface
(of the right/left-handed trefoil (Type VII) or of figure eight (Type VIII)),
see \cite{Ba, Ba3} and \cite{Y1}.
\item[(3)] {\it Sporadic examples (a), (b), (c) and (d)} \
(Type IX, X, XI and XII, respectively)
\end{enumerate}
Their surgery coefficients are also decided.
Thus we call them
{\it Berge's knots} of lens space surgery,
or {\it Berge's lens space surgeries}.
The numbering VII, $\cdots$ and XII are also used in
the recent works by Baker in \cite{Ba2, Ba3}.
It is conjectured by Gordon \cite{go1, go2} that
every knot of lens space surgery is a doubly-primitive knot.
Berge has claimed that his list of doubly-primitive knots
is complete (i.e., any doubly-primitive knot belongs to
(1), (2) or (3)), but it has not appeared.
\par
In the present paper, we are concerned with the family (1).
Its subfamily Type I consists of torus knots.
Type II consists of $2$-cables of torus knots.
Their presentations as A'Campo's divide knots
are already studied in \cite{GHY} and \cite{Y2}.
Thus our targets are Type III, IV, V and VI.
\par \medskip
\noindent
{\bf Notation.} Throughout the paper, we let
the symbol ${\mathcal X}$ denote one of these Types, i.e.,
${\mathcal X} =$ III, IV, V or VI.
\par \medskip
\noindent
To describe the knots in each Type $\mathcal{X}$, in \cite{Bg2},
Berge defined five parameters
$\delta, \varepsilon \in \{ \pm 1\}$ and $A, B, b \in {\bf Z}$
(They satisfy some certain conditions depending on $\mathcal{X}$).
We introduce two new parameters $k, t$ such that
$B, b$ are uniquely calculated from $k, t$ and vice-versa.
By $K_{\mathcal{X}}(\delta, \varepsilon, A, k, t)$, we mean
the knot defined by the parameters in Type $\mathcal{X}$.
(Type VI is slightly different from the others.)
Taking opposite $\delta$ corresponds to the mirror image of the knot.
Note that, if a lens space surgery $(K, r)$ belongs to Type ${\mathcal X}$,
$(K!, -r)$ is also a lens space surgery and
belongs to the same Type ${\mathcal X}$,
where $K!$ is the mirror image of $K$.
See Section~\ref{sec:Berge} for details on the parameters.
\par
\begin{figure}[h]
\begin{center}
\includegraphics[scale=0.6]{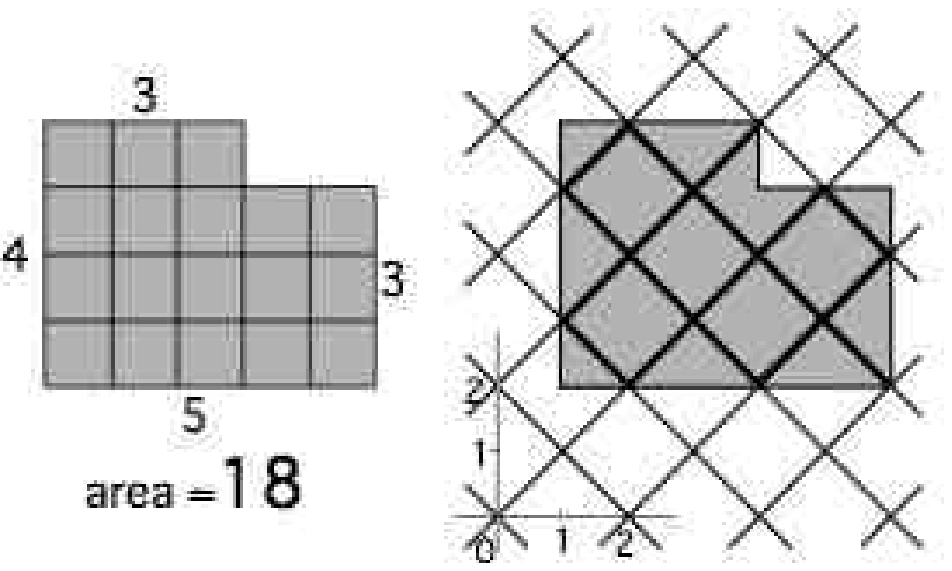}
\includegraphics[scale=0.6]{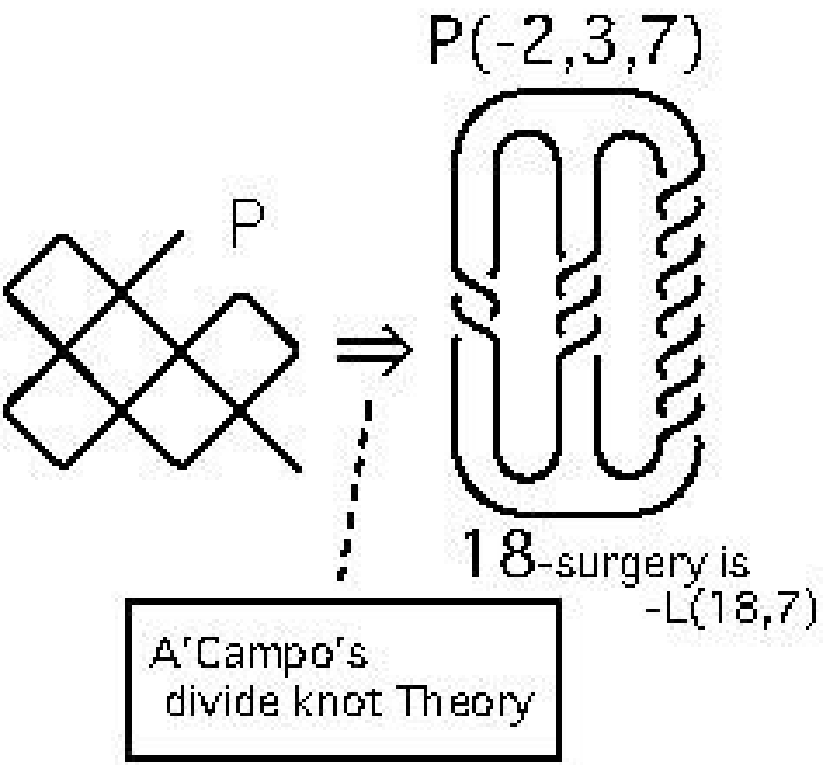}
\caption{Pretzel knot $(-2,3,7)$ with coefficient $18$}
\label{fig:237}
\end{center}
\end{figure}
The theory of A'Campo's {\it divide knots and links}
comes from singularity theory of complex curves.
The {\it divide} is (originally)
a relative, generic immersion of a 1-manifold
in a unit disk in ${\bf R}^2$.
A'Campo \cite{A1, A2, A3, A4}
formulated the way to associate to each divide $C$ a link
$L(C)$ in $S^3$.
In the present paper, we regard a PL (piecewise linear) plane curve as a divide
by smoothing the corners.
The class of divide links properly contains the class
of the links arising from isolated singularities of complex
curves, i.e., positive torus knots, and iterated torus knots
satisfying certain inequalities in their parameters.
%
%
\begin{defn}~\label{def:L}
{\rm 
Let $X$ be the $\pi /4$-lattice defined by $\{ (x,y) \vert \cos \pi x = \cos \pi y \}$
in $xy$-plane (${\bf R}^2$).
By an {\it L-shaped region}, we mean a union of two rectangles
sharing a corner and overlapping along an edge of one,
where rectangles are assumed that all edges
are parallel to either $x$-axis or $y$-axis,
and that all vertices are at lattice points ($\in {\bf Z}^2$).
We call a plane curve an {\it L-shaped curve} if
the curve $P$ is obtained as an intersection $X\cap \mathcal{L}$
of $X$ and an L-shaped region $\mathcal{L}$.
We define {\it $\textrm{area}(P)$} of an L-shaped curve $P = X\cap \mathcal{L}$
as the area (2-dim. volume) of the L-shaped region $\mathcal{L}$ defining $P$.
}
\end{defn}
%
%
See Figure~\ref{fig:237}. It is the starting example of our results.
The L-shaped curve $P = X\cap \mathcal{L}$, as a divide,
presents the pretzel knot of type $(-2,3,7)$.
Its $18$ surgery is a lens space,
which is one of the examples in \cite{FS}.
Note that the area of $P$ is equal to $18$,
the coefficient of the lens space surgery.
Our main result is:
%
%
\begin{thm}~\label{thm:main0} 
{\sl
Up to mirror image,
every Berge's knot of lens space surgery in Type III, IV, V and VI
is one of A'Campo's divide knots, and can be presented by an L-shaped curve.
}
\end{thm}
%
%

In fact, for the given parameters $\mathcal{X}$ and $(\delta , \varepsilon, A, k, t)$,
we will construct an L-shaped curve
$P_{\mathcal{X}}(\varepsilon, A, k, t)
= X \cap \mathcal{L}_{\mathcal{X}}(\varepsilon, A, k, t)$,
see Demonstration in Subsection~\ref{sbsec:demo}.
Note that opposite $\delta$ corresponds to the mirror image.

%
%
\begin{thm}~\label{thm:main1} 
{\sl
Our L-shaped curve $P_{\mathcal{X}}(\varepsilon, A, k, t)$
presents the Berge's knot
$K_{\mathcal{X}}(\delta , \varepsilon, A, k, t)$
in Type $\mathcal{X}$ ($\mathcal{X} =$ III, IV, V or VI),
up to  mirror image:
\[
L(\, P_{\mathcal{X}}(\varepsilon, A, k, t) \, )
= 
K_{\mathcal{X}}(1, \varepsilon, A, k, t),
\textrm{ or its mirror image }
K_{\mathcal{X}}(-1, \varepsilon, A, k, t).
\]
\par
One of $K_{\mathcal{X}}(\pm 1, \varepsilon, A, k, t)$ is
presented by a positive braid (say $w$) and the other
is by a negative one (the inverse $w^{-1}$).
The divide knot $L(\, P_{\mathcal{X}}(\varepsilon, A, k, t) \, )$ is
exactly equal to the positive one, but
the choice at $\delta$ ($1$ or $-1$) depends on
$\mathcal{X}$ and $(\varepsilon, A, k, t)$.
}
\end{thm}
%
%

Next, we study the surgery coefficients.
By the Cyclic Surgery Theorem of
Culler--Gordon--Luecke--Shalen \cite{CGLS},
if a hyperbolic knot $K$ admits a lens space surgery,
then the coefficient is integral.
By $\textrm{coef}( K_{\mathcal{X}}(\delta , \varepsilon, A, k, t) )$,
we denote the surgery coefficient of the lens space surgery
of the knot \underline{as in Type $\mathcal{X}$}.
Note that there exist some hyperbolic knots
that have two coefficients of lens space surgery
(such coefficients are proved to be consecutive in \cite{CGLS}),
and belong to different Types as the pairs with the coefficients.
It is the reason why we state
\lq\lq as in Type $\mathcal{X}$\rq\rq.
%
%
\begin{thm}~\label{thm:main2} 
{\sl
Under the correspondence in Theorem~\ref{thm:main1},
the area of the L-shaped curve $P$ is equal to
(the absolute value of) the coefficient of
the lens space surgery of $L(P)$ as in Type $\mathcal{X}$,
or is greater by one:
\[
\textrm{area}(P) - \vert \textrm{coef}(L(P)) \vert = 0 \textrm{ or } 1.
\]
}
\end{thm}
%
%
This theorem will be proved as Lemma~\ref{lem:v-c},
in which we will decide the choice ($0$ or $1$)
by the parameters.
We will prove that $\textrm{coef}(L(P)) > 0$ in Lemma~\ref{lem:pos}.
Thus we will change $\vert \textrm{coef}(L(P)) \vert$
to $\textrm{coef}(L(P))$ in Lemma~\ref{lem:v-c}.

Theorem~\ref{thm:main0} can be proved
by combination of Lemma~\ref{lem:Berge} and
Lemma~\ref{lem:mirror}.
But the aim of the present paper is the construction of
the L-shaped curves 
using the operation {\it adding squares} on L-shaped curves
in Section~\ref{sec:L},
and studying the knots and the family of knots by them.
Some knots are
obtained from other knots (possibly in other Types)
by some twistings.
By our method adding squares,
we can search such pairs, and check such relations easily,
see Section~\ref{sec:Obs}.

\medskip
Here we survey on divide presentation of
the other Berge's knots (in Type VII and later), shortly.
All knots are considered up to mirror image.
Type VII consists of the knots
that the author \cite{Y1} gave
L-shaped curve presentations first.
Type VIII contains some knots that is hard
(in the author's opinion) to decide
whether it is a divide knot or not,
and (if it is) to present by a divide.
The author has shown that
every sporadic knot (in Type IX and later)
is a divide knot and has shown a method to construct the divide.
But he does not know whether it can be presented by an
\underline{L-shaped} curve or not.

Note that there exist a family of
L-shaped divide knot whose $\textrm{area}(P)$-surgery is not
a lens space \cite{Y3}, but such L-shaped divide knots tend to have
exceptional Dehn surgeries, to the author's knowledge \cite{Y2, Y4}.
\par

\medskip
This paper is organized as follows: In the next section, we review
Berge's knots and their parameters in detail.
In Section~\ref{sec:ACampo},
we review A'Campo's divide knot theory and define L-shaped curves.
In Section~\ref{sec:L},
developing a method {\it adding squares},
we construct L-shaped curves (regions) for Berge's knots.
In Section~\ref{sec:Proof},
we will prove Theorem~\ref{thm:main1} and
Lemma~\ref{lem:v-c}, the precise version of Theorem~\ref{thm:main2}.
Finally, in Section~\ref{sec:Obs}, we study some applications,
advantages to present Berge's knots as divide knots.
We place Tables~\ref{tbl:Para} and \ref{tbl:Lregion}
after 
the reference list for the reader's convenience. 
\par

\section{Berge's knots of Type III, IV, V and VI}~\label{sec:Berge}
We recall the Berge's parametrization of knots in Type ${\mathcal X}$.
We use his original parameters
$\delta, \varepsilon, A, B, b$ and a constant $a$
($:= 0$ or $1$) defined in \cite{Bg2},
and introduce two new parameters $k$ and $t$.
\par
We start with the following:
\begin{enumerate}
\item[(1)]
$\cdot$ $\delta$ and $\varepsilon$ are signs ($\in \{ \pm 1\}$).
The opposite $\delta$ corresponds to the mirror image. 
\par
$\cdot$ $A$ is a positive integer, whose range and parity (even or odd)
depends on ${\mathcal X}$, 
\par
$\cdot$ $k$ runs in ${\bf N}_{\geq0}$.
$B$ is decided by $(\varepsilon, A, k)$
and satisfies that $0 < 2A \leq B$.
\par
$\cdot$ $b, t \in {\bf Z}$. They can be negative.
\item[(2)]
The parameters $k, t $ are introduced instead of
the conditions in \cite{Bg2} written by sentences and
by congruences, respectively.
For example, instead of
\lq\lq $(B + \varepsilon)/A$ is an odd integer\rq\rq \ %
in \cite[Table 3(p.15)]{Bg2},
we set $B = A(3+2k) -\varepsilon$.
Instead of \lq\lq $b \equiv -2 \varepsilon \delta A$ (mod $B$)\rq\rq \ %
in \cite[Table 3(p.15)]{Bg2},
we set $b = - \delta \varepsilon (2A + tB)$.
These are the relations between $(B, b)$ and $(k, t)$ in Type III.
They are similar in other Types, but slightly different, see Table~\ref{tbl:Para}(1).
\item[(3)] 
The independent parameters are
$(\delta, \varepsilon, A, k, t)$  in Type III, IV and V,
but is $(\delta, A, t)$  in Type VI.
We formally regard the latter as $(\delta, \varepsilon, A, k, t) = (\delta, -1, A, 0, t)$,
i.e., we fix $\varepsilon := -1, k := 0$ in Type VI, for the convenience.
\end{enumerate}
\noindent
{\bf Notation.}
By $K_{{\mathcal X}}(\delta, \varepsilon, A, k, t)$,
we denote the knot parametrized
as $(\delta, \varepsilon, A, k, t)$ in Type ${\mathcal X}$,
by Berge in \cite{Bg2}.
\par
\medskip
Now we go into the detail. See Table~\ref{tbl:Para}(1), (2) and (3).
In Table~\ref{tbl:Para}(1),
we define $B$ and $b$ (depending on ${\mathcal X}$), using
a temporary parameter $l$.
For fixed $A$ and $\varepsilon$, the possible values of $l$ are
in an arithmetic sequence, depending on $\mathcal{X}$.
We parametrize the sequence by $k \in {\bf N}_{\geq 0}$ 
as in Table~\ref{tbl:Para}(2).
In every case, the surgery coefficient is $bB + \delta A$,
where $B$ depends on $\varepsilon, A, k$.
In Table~\ref{tbl:Para}(3),
we deform the coefficients into the form including the terms
$+ kA^2 + tB^2$ (or $+ k(2A)^2 + tB^2$ in Type III).
These are related to our method adding squares in
Section~\ref{sec:L}.
Note that, if a knot $K$ with coefficient $r$ belongs to Type ${\mathcal X}$,
its mirror image $K!$ with $-r$ (i.e., opposite $\delta$)
also belongs to the same Type ${\mathcal X}$.
\par
\begin{figure}[h]
\begin{center}
\includegraphics[scale=0.6]{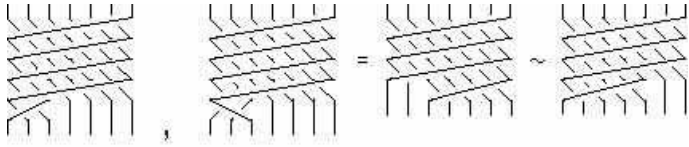}
\caption{
$W(7)^4 W(3)$ and $W(7)^4 W(3)^{-1} \sim W(7)^3 W(5) $
}
\label{fig:Wn}
\end{center}
\end{figure}

Using these parameters, in \cite{Bg2},
Berge has already given the braid presentations
of these knots:
%
%
\begin{lem}~\label{lem:Berge} {\rm (Berge \cite{Bg2})} \quad 
{\sl
Every knot $K_{{\mathcal X}}(\delta, \varepsilon, A, k, t)$
is presented as the closure of the braid
$W(B)^{b} W(A+1-a)^{\delta}$ of index $B$,
where $W(n) = \sigma_{n-1} \sigma_{n-2} \cdots  \sigma_{1}$,
see Figure~\ref{fig:Wn}.
}
\end{lem}
%
%

\begin{figure}[h]
\begin{center}
\includegraphics[scale=0.5]{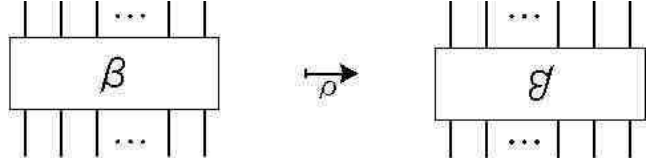}
\caption{$\pi$-rotation $\rho$}
\label{fig:Pirot}
\end{center}
\end{figure}

%
%
\begin{defn}~\label{def:rot}
{\rm 
We define an anti-homomorphic (i.e.,
$\rho(\beta_1 \beta_2) = \rho(\beta_2) \rho (\beta_1)$)
involution {\it $\pi$-rotation}
$\rho$ on the braid group of index $n$ by
extending $\rho (\sigma_i) = \sigma_{n-i}$, see
Figure~\ref{fig:Pirot}.
By $\beta' \stackrel{\rho}{=} \beta$, we mean
$\beta' = \rho (\beta)$ and equivalently $\rho (\beta') = \beta$.

In this opportunity, we define another notation
$\beta' \sim \beta$ as that
the closure of $\beta'$ is the same
knot or link to that of $\beta$.
Note that $\beta' \stackrel{\rho}{=} \beta$ up to conjugate implies
$\beta' \sim \beta$.
}
\end{defn}
%
%
Our L-shaped divide knots are always presented by positive braids,
see Section~\ref{sec:ACampo}, while any divide knot is a closure of
a strongly quai-positive braid (Lemma~\ref{lem:divides}(7)).
Thus, first, if $b<0$, we take the mirror image
(i.e., change the sign $\delta$, then $b$ becomes to $-b >0$) and next,
we use the following lemma if it is necessary.
%
%
\begin{lem}~\label{lem:Wb-1}
{\sl
Let $a_1, a_2$ and $c$ be positive integers with $a_1 < a_2$.
The closure of the braid $W(a_2)^c W(a_1)^{-1}$ of index $a_2$
is the same knot to that of $W(a_2)^{c-1} W(a_2 - a_1 +1)$.
}
\end{lem}
%
%
{\it Proof.}
\begin{eqnarray*}
W(a_2) W(a_1)^{ -1}
& = & 
\sigma_{a_2-1} \sigma_{a_2-2} \cdots \sigma_{a_1}
\cdots  \sigma_2 \sigma_1
(\sigma_{a_1-1} 
\cdots \sigma_2 \sigma_1)^{-1} \\
& = & 
\sigma_{a_2-1} \sigma_{a_2-2} \cdots \sigma_{a_1} \\
& \stackrel{\rho}{=} & 
W(a_2 - a_1+1).
\end{eqnarray*}
Since $\pi$-rotation $\rho$ is anti-homomorphic, and 
the braid $W(n)$ of index $n$ ($=a_2$) is fixed by $\rho$, 
we have
\begin{eqnarray*}
W(a_2)^c W(a_1)^{-1}
& = &
W(a_2)^{c-1} \cdot W(a_2) W(a_1)^{-1} \\
& \stackrel{\rho}{=} & 
W(a_2 - a_1+1) W(a_2)^{c-1},
\end{eqnarray*}
which is conjugate to $W(a_2)^{c-1} W(a_2 - a_1 +1)$.
\qed
\par \medskip
Figure~\ref{fig:Wn} illustrates Lemma~\ref{lem:Wb-1}.
$K_{\textrm{III}}(1,-1, 2, 0, 0)$ is
presented by the braids $W(7)^4 W(3)$,
and $K_{\textrm{III}}(-1,1, 2, 0, 0)$ is presented by
$W(7)^4 W(3)^{-1} \sim W(7)^3 W(5)$. 
\par

\section{L-shaped curves and A'Campo's divide knots}~\label{sec:ACampo}
The theory of A'Campo's {\it divide knots and links} \cite{A1, A2, A3, A4}
comes from singularity theory of complex curves.
It is a method to associate to each divide (a plane curve) $C$
a link $L(C)$ in the $3$-dimensional sphere $S^3$.
The original definition of divide knots in \cite{A1} is
differential-geometric.
Hirasawa \cite{Hi} visualized the construction.
%
We are concerned with the plane curves of special type, called
\lq\lq L-shaped curves\rq\rq,
see Subsection~\ref{sbsec:L}.
For such special curves,
we can use another method introduced by
Couture--Perron \cite{CP},
see Subsection~\ref{sbsec:Lknot}.

\subsection{A'Campo's divide knots}~\label{sbsec:ACknot}
\begin{figure}[htbp]
\begin{center}
\includegraphics[scale=0.6]{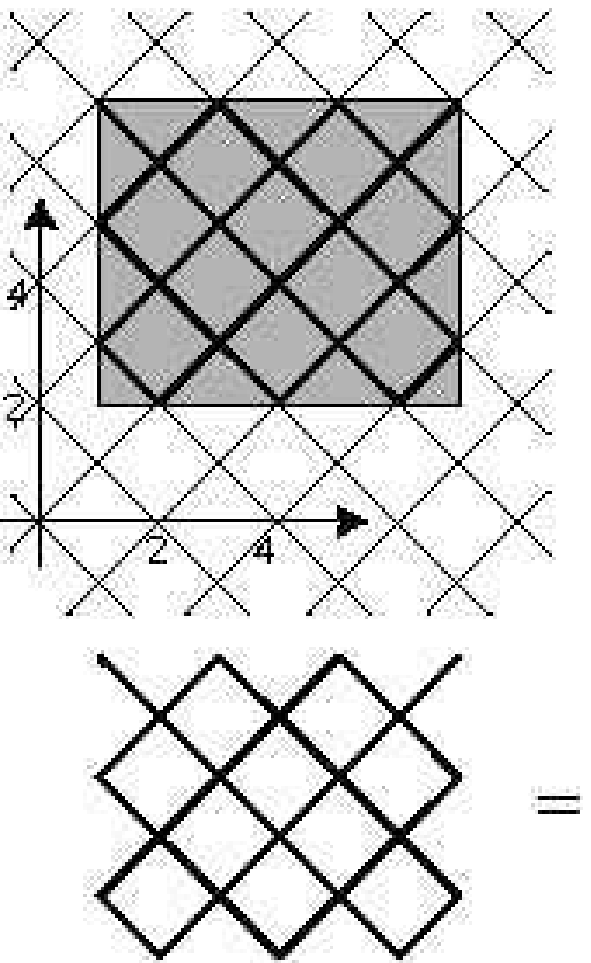}
\includegraphics[scale=0.6]{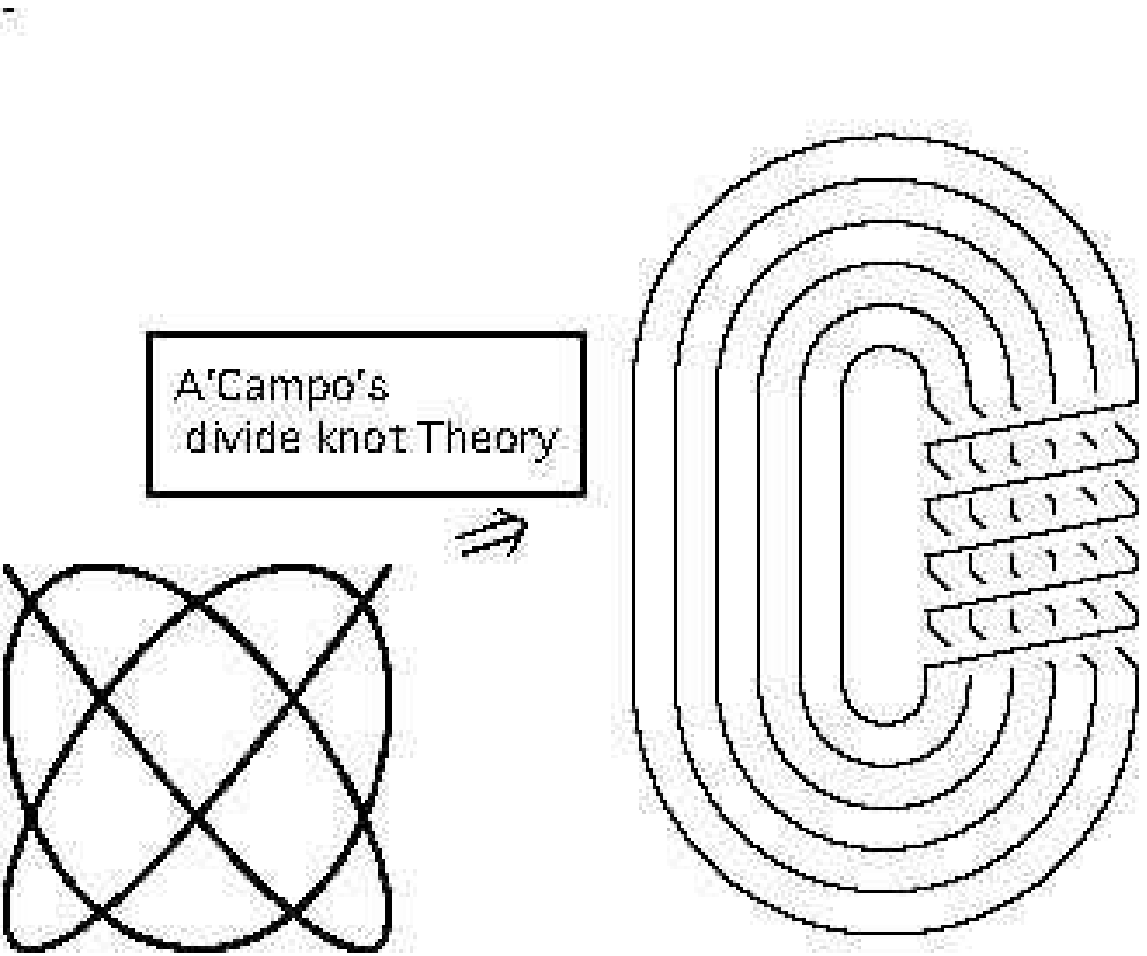}
\caption{A billiard curve presents a torus knot (ex. $T(6,5)$)}
\label{fig:Torus}
\end{center}
\end{figure}

\begin{figure}[h]
\begin{center}
\includegraphics[scale=0.6]{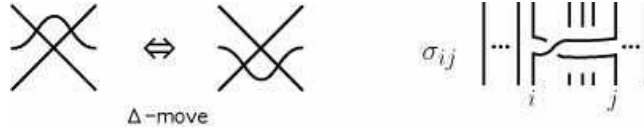} \par
\caption{Basics on divide knots}
\label{fig:div}
\end{center}
\end{figure}

We start with the typical example of divide knots,
see Figure~\ref{fig:Torus}:
%
%
\begin{lem}~\label{lem:GHY}
{\rm (Goda--Hirasawa--Y \cite{GHY}, see also \cite{AGV, Gu})} \quad 
{\sl 
Let $a,b$ be a pair of positive integers
and $\mathcal{R}(a,b)$ be an $a \times b$-rectangle region.
A plane curve defined by $X \cap \mathcal{R}(a,b)$
(a billiard curve of type $B(a,b)$)
presents the torus link $T(a,b)$ as a divide.
}
\end{lem}
%
%
Some characterizations of (general) divide knots and links are known,
and some topological invariants $L(P)$
can be gotten from the divide $P$ directly.
Here, we list some of them.
%
%
\begin{lem}~\label{lem:divides}
{\rm  ((1)--(6) by A'Campo \cite{A2},
(7) by Hirasawa \cite{Hi}, Rudolph \cite{R})}
{\sl
\begin{enumerate}
\item[(1)]
$L(P)$ is a knot (i.e., connected) if and if only
$P$ is an immersed arc.
\item[(2)]
If $L(P)$ is a knot, the unknotting number, genus and $4$-genus
of $L(P)$ are all equal to
the number $d(P)$ of the double points of $P$.
\item[(3)]
If $P = P_1 \cup P_2$ is the image of an immersion of two arcs,
then the linking number of the two component link
$L(P) = L(P_1) \cup L(P_2)$ is equal to
the number of the intersection points between $P_1$ and $P_2$.
\item[(4)]
If $P$ is connected, then $L(P)$ is {\it fibered}.
\item[(5)]
A divide $P$ and its mirror image $P!$ present
the same knot or link: $L(P!) =L(P)$.
\item[(6)]
If $P_1$ and $P_2$ are related by some $\Delta$-moves,
then the links $L(P_1)$ and $L(P_2)$ are isotopic:
If $P_1 \sim_{\Delta} P_2$ then $L(P_1) = L(P_2)$,
see Figure~\ref{fig:div}.
\item[(7)]
Any divide knot is a closure of a
{\it strongly quasi-positive} braid,  i.e.,
a product of some $\sigma_{ij}$ in Figure~\ref{fig:div}.
\end{enumerate}
}
\end{lem}
%
%
For theory of divide knots, see also Rudolph's \lq\lq {\bf C}-link\rq\rq \ %
in \cite{R} and \cite{Ch, HW}.

\subsection{Preliminary on L-shaped curves}~\label{sbsec:L}
First, we parametrize L-shaped regions
by four positive integers $a_1, a_2, b_1, b_2$ that satisfy
$a_1 < a_2$ and $b_1<b_2$, see Figure~\ref{fig:Lreg}: 
\par
%
%
\begin{defn}~\label{def:L-sh}
{\rm (L-shaped region at the origin) \quad 
In $xy$-plane, we define
\begin{eqnarray*}
L[a_1, a_2; b_1, b_2]
& := &
\{
(x,y)
\vert \,
0 \leq x \leq a_1,
0 \leq y \leq b_2
\} \\
 & &
\cup
\{
(x,y)
\vert \,
0 \leq x \leq a_2,
0 \leq y \leq b_1
\}.
\end{eqnarray*}
\begin{figure}[h]
\begin{center}
\includegraphics[scale=0.6]{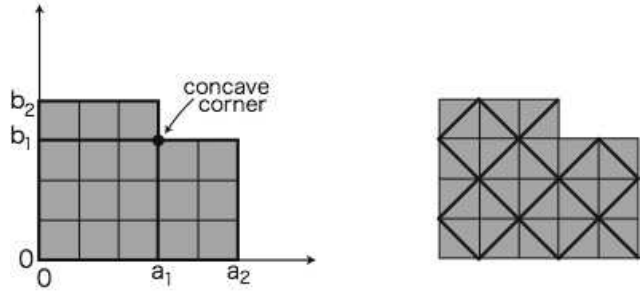}
\caption{L-shaped region $L[3,5; 3,4]$}
\label{fig:Lreg}
\end{center}
\end{figure}
By {\it concave corner},
we mean the point (of the region) at the coordinate
$(a_1,b_1)$ in the definition above.
We will call not only $L :=L[a_1, a_2; b_1, b_2]$
but also its transformations
${\mathcal L} = T(L)$ an
{\it L-shaped region of type $[a_1, a_2; b_1, b_2]$},
where $T$ is a transformation in $xy$-plane generated by
the reflection $r_X$ along the $x$-axis (Lemma~\ref{lem:divides}(5)),
the rotation $R$ by $\pi /2$ about a lattice point and
the parallel translation $+ \vec{n}$ by
a lattice point $\vec{n}$ ($\in {\bf Z}^2$).
}
\end{defn}
%
%

Let $X$ be the $\pi /4$-lattice
defined by $\{ (x,y) \vert \cos \pi x = \cos \pi y \}$
in $xy$-plane.
A lattice point $(m,n)$ ($\in {\bf Z}^2$) is called
{\it even} (or {\it odd}, resp.)
if $m+n$ is even (or odd).
We are concerned only with the case that
the intersection $X \cap \mathcal{L}$ is
the image of a generic immersed arc.
Thus, we always control $+\vec{n}$ and assume that
\begin{center}
($\ast$) The concave point of an L-shaped region ${\mathcal L}$ is
placed at an odd point.
\end{center}
Assuming ($\ast$), the parameter $[a_1, a_2; b_1, b_2]$
defines a unique plane curve up to isotopy, i.e.,
it depends on neither $r_X, R$ nor
translations keeping even/odd points.
We call the corresponding plane curve an
{\it L-shaped curve of type $[a_1, a_2; b_1, b_2]$}.
Of course, for an L-shaped curve $P$ of type $[a_1, a_2; b_1, b_2]$,
we have
\[
\textrm{area}(P) = a_2 b_1 + a_1 b_2 - a_1 b_1.
\]
On the other hand, the number $d(L)$ of double points of $L$ is
\[
d(P) = \{ a_2 (b_1-1) + b_2 (a_1 -1) - a_1 b_1 +1 \}/2,
\]
because double points
are the even points of the interior of the L-shaped region.
\par
The condition ($\ast$) is not sufficient for $X \cap \mathcal{L}$ to be
the image of an immersed arc.
In fact, it possibly consists of multiple components or
contains some circle components.
\par
\medskip
The following proposition follows from
Lemma~\ref{lem:divides}(2).
%
%
\begin{prop}~\label{prop:GT}
{\sl If the L-shaped curve
$P =X \cap \mathcal{L}$ of type $[a_1, a_2; b_1, b_2]$
with the assumption ($\ast$) is an immersed arc, then
the genus $g(L(P))$ of the divide knot $L(P)$ is
(the unknotting number, and the $4$-genus are also)
equal to the number $d(P)$ of the double points of $P$:
\[
g(L(P)) = \{ a_2 (b_1-1) + b_2 (a_1 -1) - a_1 b_1 +1 \}/2.
\]
Thus, it holds that $\textrm{area}(P) - 2g(L(P)) = a_2 +b_2 -1$.}
\end{prop}
%
%
\par

\subsection{L-shaped divide knots}~\label{sbsec:Lknot}
In \cite{CP} Couture and Perron pointed out a method
to get the braid presentation from the divide (the plane curve)
in the restricted cases, called \lq\lq ordered Morse\rq\rq \ %
divides.
Our L-shaped curves are all ordered Morse,
thus we can apply their method.
It is a special case of Hirasawa's method in \cite{Hi}.
%
%
\begin{lem}~\label{lem:QtwistT}
{\sl
The divide link presented by the L-shaped curve of type
$[a_1, a_2; b_1, b_2]$ is the closure of
the braid $W(a_2)^{b_1}  W(a_1)^{b_2 - b_1}$ of index $a_2$,
where $W(n) = \sigma_{n-1} \sigma_{n-2} \cdots  \sigma_{1}$.
}
\end{lem}
%
%
Such a link should be regarded as a \lq\lq rationally twisted\rq\rq \ %
torus link in the following sense:
The link is obtained by a \lq\lq $(b_2-b_1)/a_1$ twist\rq\rq \ %
of the parallel $a_1$ strings in $a_2$ strings
of torus link $T(a_2, b_1)$ in the standard position $W(a_2)^{b_1}$.
%
%
\begin{ex}~\label{ex:237}
{\rm
The divide knot presented by the L-shaped curve $[3,5;3,4]$
is the closure of the braid
$(\sigma_2 \sigma_4 \sigma_1 \sigma_3)^3 \sigma_2 \sigma_1$
(conjugate to $(\sigma_4 \sigma_3 \sigma_2 \sigma_1)^3 \sigma_2 \sigma_1$)
of index $5$, which is $P(-2,3,7)$.
}
\end{ex}
%
%
{\it Proof.}
\quad 
First, we define the words $o(n)$ and $e(n)$ in the braid group
of index $a_2$ as follows:
\[
e(n) := \prod_{i \textrm{: even}, \, i <n} \sigma_i \ , \qquad
o(n) := \prod_{i \textrm{: odd}, \, i <n} \sigma_i \ , \qquad
\]
where $n$ is a positive integer less than or equal to the index.
Note that $\sigma_i$ and $\sigma_j$ are commutative
if $i$ and $j$ have the same parity.
If $j \leq k < l$ and $k$ is even, $o(k)^{-1}o(l)$
is a product of $\sigma_i$'s with $i \geq k+1$, thus
is commutative with both $o(j)$ and $e(j)$.
Similarly, if $j \leq k < l$ and $k$ is odd, then
$e(k)^{-1}e(l)$ is commutative with both $o(j)$ and $e(j)$.

\begin{figure}[h]
\begin{center}
\includegraphics[scale=0.6]{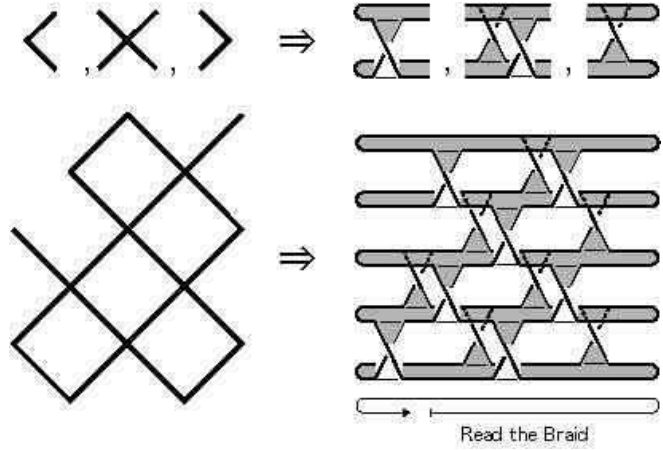}
\caption{Couture-Perron's method in the case $[3, 5; 3, 4]$ ($\pi/2$ rotated)}
\label{fig:CP}
\end{center}
\end{figure}

In the case of L-shaped curves, Couture--Perron's method is
summarized as the algorithm in Figure~\ref{fig:CP}.
By direct application of the algorithm to
the L-shaped curve of type, we have:
\par
\begin{figure}[h]
\begin{center}
\includegraphics[scale=0.6]{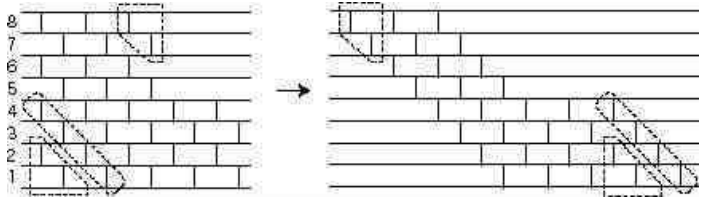}
\caption{Braid of L-shaped curve $[5, 9; 3, 5]$}
\label{fig:Amida}
\end{center}
\end{figure}

\noindent
{\bf Claim~1} {\rm \cite{CP}}
\quad 
{\sl 
The L-shaped curve of type $[a_1, a_2; b_1, b_2]$
presents the closure of the braid
\[
B[a_1, a_2; b_1, b_2] := \begin{cases}
\begin{pmatrix} e(a_2)o(a_2) \end{pmatrix}^{b_1}
\begin{pmatrix} e(a_1)o(a_1) \end{pmatrix}^{b_2 -b_1}
& \textrm{if $a_1$ is odd, }\\
\begin{pmatrix} o(a_2)e(a_2) \end{pmatrix}^{b_1}
\begin{pmatrix} o(a_1)e(a_1) \end{pmatrix}^{b_2 -b_1}
& \textrm{if $a_1$ is even.}
\end{cases}
\]
}
The key idea of the rest of the proof is in Figure~\ref{fig:Amida}.
Let $G(n)$ be as follows.
\[
G(n) := \begin{cases}
e(n+1) o(n) e(n-1) o(n-2) e(n-3) \cdots e(4)o(3)
& \textrm{if $n$ is odd, } \\
o(n+1) e(n) o(n-1) e(n-2) o(n-3) \cdots e(4)o(3)
& \textrm{ if $n$ is even.}
\end{cases}
\]
It is a product of some $\sigma_i$'s with $i < n$. 
\par \medskip
\noindent
{\bf Claim~2}
{\sl
\begin{eqnarray*}
G(n-2)^{-1} e(n)o(n) \ G(n-2) = W(n)
& \textrm{ if $n$ is odd,} \\
G(n-2)^{-1} o(n)e(n) \ G(n-2) = W(n)
& \textrm{ if $n$ is even.}
\end{eqnarray*}
}
In fact, if $n$ is odd,
$e(n-1)^{-1}e(n) = \sigma_{n-1}$ and it commutes
with $o(n-2)^{-1}$.
Next, $o(n-2)^{-1} o(n) = \sigma_{n-2}$
and it commutes with $e(n-3)^{-1}$.
We repeat such reductions
until $o(3)^{-1} o(5) = \sigma_3$ inductively,
and end with $e(4)o(3) = \sigma_2 \sigma_1$.
The other case is proved similarly. 
\par
Next, we set $H(a_1, a_2) := \rho( G(a_2 - a_1 -1) )$,
where $\rho$ is $\pi$-rotation in Definition~\ref{def:rot}.
Then, $H(a_1, a_2)$ is a product of some $\sigma_i$'s with $i > a_1+1$.
Thus we have
\par \medskip
\noindent
{\bf Claim~3} \
{\sl
$H(a_1, a_2)$ commutes with $e(a_1), o(a_1)$ and $G(a_1-2)$.
}
\par \medskip
Let $\Omega(a_1, a_2) := H(a_1, a_2)^{-1}G(a_1-2)$.
By Claim~2 and 3, It holds 
\par \medskip
\noindent
{\bf Claim~4}
{\sl
\begin{eqnarray*}
\Omega(a_1, a_2)^{-1}
e(a_1) o(a_1)
\Omega(a_1, a_2) = W(a_1)
& \textrm{ if $a_1$ is odd,} \\
\Omega(a_1, a_2)^{-1}
o(a_1) e(a_1)
\Omega(a_1, a_2) = W(a_1)
& \textrm{ if $a_1$ is even.}
\end{eqnarray*}
}
The following is the most troublesome step. 
\par \medskip
\noindent
{\bf Claim~5}
{\sl
\begin{eqnarray*}
\Omega(a_1, a_2)^{-1} e(a_2)o(a_2) \ 
\Omega(a_1, a_2) = W(a_1)W(a_2)W(a_1)^{-1}
& \textrm{ if $a_1$ is odd,} \\
\Omega(a_1, a_2)^{-1} o(a_2)e(a_2) \ 
\Omega(a_1, a_2) = W(a_1)W(a_2)W(a_1)^{-1}
& \textrm{ if $a_1$ is even.}
\end{eqnarray*}
}
To prove Claim~5, we divide the braid of index $a_2$ into two parts,
lower and higher parts along the $a_1$-th string.
Here we denote $a_2-a_1 +1$ by $\overline{a_1}$.
In the case $a_1$ is odd,
\[
e(a_2) o(a_2) =
\begin{cases}
e(a_1) o(a_1) \,
\rho (
e(\overline{a_1}) o(\overline{a_1}) )
& \textrm{ if $a_2$ is odd,}
\\
e(a_1) o(a_1) \,
\rho (
o(\overline{a_1}) e(\overline{a_1}) )
& \textrm{ if $a_2$ is even.}
\end{cases}
\]
In the former case, by Claim~3,
the conjugation of $e(a_2) o(a_2)$ by
$\Omega(a_1, a_2)$ is divided as
the product of
that of $e(a_1) o(a_1)$ by $G(a_1-2)$
and
that of $\rho (
e(\overline{a_1}) o(\overline{a_1}) )$
by $H(a_1, a_2)^{-1}$.
Since $H(a_1, a_2)$ is defined as $\rho( G(\overline{a_1} -2) )$,
it holds
\begin{eqnarray*}
H(a_1, a_2) \,
\rho (
e(\overline{a_1}) o(\overline{a_1}) ) \,
H(a_1, a_2)^{-1}
& = &
\rho (
G(\overline{a_1} -2)^{-1}
e(\overline{a_1}) o(\overline{a_1})
G(\overline{a_1} -2)
) \\
& = &
\rho (
W(\overline{a_1})
)
\\
& = &
W(a_2) W(a_1)^{-1}
\end{eqnarray*}
Here, we use Claim~2 with odd $n = \overline{a_1}$.
In the other cases, including $a_2$ is even,
the proofs are similar. 
\par
Finally, by Claim~1, 4 and 5,
\begin{eqnarray*}
\Omega(a_1, a_2)^{-1} B[a_1, a_2; b_1, b_2] \ \Omega(a_1, a_2)
& = & (W(a_1)W(a_2)W(a_1)^{-1})^{b_1} W(a_1)^{b_2-b_1} \\
& = & W(a_1)W(a_2)^{b_1}W(a_1)^{-1} W(a_1)^{b_2-b_1}
\end{eqnarray*}
\par
It is conjugate to $W(a_2)^{b_1} W(a_1)^{b_2-b_1}$.
The proof of Lemma~\ref{lem:QtwistT} is complete.
\qed
\par \medskip

By the symmetry between the L-shaped curve
of type $[a_1, a_2; b_1, b_2]$ and that of type $[b_1, b_2; a_1, a_2]$,
we have an extension of the well-known
symmetry $T(b,a) = T(a,b)$ of torus knots.
%
%
\begin{cor}~\label{cor:sym}
{\sl
The closures of the braids
\[
W(b_2)^{a_1}  W(b_1)^{a_2 - a_1}
\textrm{ of index $b_2$}
\quad \textrm{and} \quad
W(a_2)^{b_1}  W(a_1)^{b_2 - b_1}
\textrm{ of index $a_2$}
\]
define the same link.
}
\end{cor}
%
%
\begin{lem}~\label{lem:mirror}
{\sl
Let $a_1, a_2$ and $c$ be positive integers with $a_1 < a_2$,
and $\delta$ be a sign ($\in \{ \pm 1\}$).
Then, the knot of the closure of the braid of type
$W(a_2)^{\pm c} W(a_1)^{\delta}$ of index $a_2$
is presented as a divide knot presented by an
L-shaped curve, up to mirror image:
\begin{enumerate}
\item[($++$)] The knot $W(a_2)^c W(a_1)$
is presented by the L-shaped curve
$[a_1, a_2; c, c+1]$.
\item[($+-$)] The knot $W(a_2)^c W(a_1)^{-1}$
is presented by the L-shaped curve
$[a_2 - a_1+1, a_2; c-1, c]$.
\item[($--$)] The knot $W(a_2)^{-c} W(a_1)^{-1}$
is the mirror image of the knot presented by the L-shaped curve
$[a_1, a_2; c, c+1]$.
\item[($-+$)] The knot $W(a_2)^{-c} W(a_1)$
is the mirror image of the knot presented by the L-shaped curve
$[a_2 - a_1+1, a_2; c-1, c]$. 
\par
\end{enumerate}
}
\end{lem}
%
%
{\it Proof.}
\quad 
The case ($++$) in the lemma follows from Lemma~\ref{lem:QtwistT}
directly, and ($+-$) follows from
Lemma~\ref{lem:Wb-1} and Lemma~\ref{lem:QtwistT}.
The cases ($--$) and ($-+$)
follow from ($++$) and ($+-$) respectively, since,
if a knot $K$ is the closure of the braid $w$, then
the mirror image $K!$ is that of the inverse $w^{-1}$, in general.
\qed
\par \medskip
By Lemma~\ref{lem:Berge} and Lemma~\ref{lem:mirror},
Theorem~\ref{thm:main0} is already proved:
Up to mirror image,
every Berge's knot of lens space surgery in Type III, IV, V and VI
is one of A'Campo's divide knots, and can be presented
by an L-shaped curve.
\par

We end this section with referring
the fiberedness of Berge's knots.
By Theorem~\ref{thm:main0} and fiberedness of divide knots in
Lemma~\ref{lem:divides}(4), we can show
%
%
\begin{cor}~\label{cor:fib}
{\rm (Teragaito \cite{Te}, Ozsv\'{a}th--Sz\'{a}bo\cite{OSz})} \quad 
{\sl Every Berge's knot of lens space surgery in Type 
III, IV, V and VI is fibered.
}
\end{cor}
%
%
This corollary can be proved by
Lemma~\ref{lem:Berge}, ~\ref{lem:Wb-1}
and the fact that knots presented by positive (or negative) braids
are fibered \cite{S}.
In fact, Teragaito \cite{Te} (see \cite[\S 5.7]{HM}),
has shown that every Berge's knots
(including Type VII, ..., XII, see Section~\ref{sec:intro}) is
fibered, by proving the braid positivity.
Ozsv\'{a}th--Sz\'{a}bo \cite[\S 5]{OSz} also proved
fiberedness of every Berge's knot 
from another view points.
\par

\section{Construction of L-shaped curves}~\label{sec:L}
We define the operation {\it adding squares} on L-shaped curves
(via L-shaped regions), its drawing notations, and
explain how to construct the L-shaped curves
$P_{\mathcal{X}}(\varepsilon, A, k, t)
= X \cap \mathcal{L}_{\mathcal{X}}(\varepsilon, A, k, t)$
for Berge's knots in Type $\mathcal{X}$,
given by the parameters $(\delta , \varepsilon, A, k, t)$.
From now on, we consider only L-shaped divide \underline{knots},
i.e., the case that the L-shaped curve $P$
(with the assumption ($\ast$) in Section~3.2)
is the image of an immersed arc (Lemma~\ref{lem:divides}(1)).
\par

\subsection{Adding squares I}~\label{sbsec:AddS}
\begin{figure}[h]
\begin{center}
\includegraphics[scale=0.4]{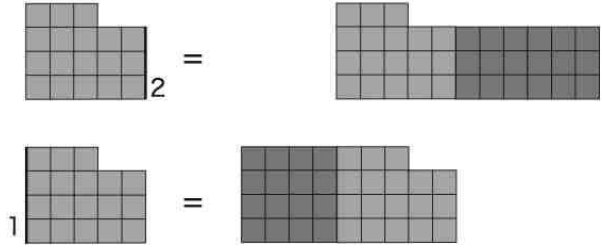}
\caption{Adding squares I}
\label{fig:LwithSq1}
\end{center}
\end{figure}

We start with the following:
%
%
\begin{defn}~\label{def:AddS1}
{\rm 
For a positive integer $n$,
we call the operation on L-shaped curves, changing
from that of type $[a_1, a_2; b_1, b_2]$
to $[a_1, a_2 +nb_1; b_1, b_2]$ or
$[a_1+ nb_2, a_2+nb_2; b_1, b_2]$ {\it adding $n$ squares},
see examples in Figure~\ref{fig:LwithSq1}.
As a drawing notation,
we specify the edge along which the squares are added,
and write $n$ near the edge.
By the symmetry in Corollary~\ref{cor:sym},
we also call the changing from $[a_1, a_2; b_1, b_2]$
to $[a_1, a_2; b_1, b_2 + n a_2]$ or
$[a_1, a_2; b_1+ na_2, b_2 + n a_2]$
adding $n$ squares.
}
\end{defn}
%
%
%
%
\begin{lem}~\label{lem:pureAS}
{\sl
Adding $n$ squares on an L-shaped curve $P$
along an edge (of the region)
corresponds to positive (i.e., right-handed) $n$ full-twists
on the divide knot $L(P)$ along the unknot defined by the edge.
}
\end{lem}
%
%
This lemma is proved by using the braid presentation
in Lemma~\ref{lem:QtwistT}.
Note that a full-twist is in the center of the braid group.
\begin{figure}[h]
\begin{center}
\includegraphics[scale=0.7]{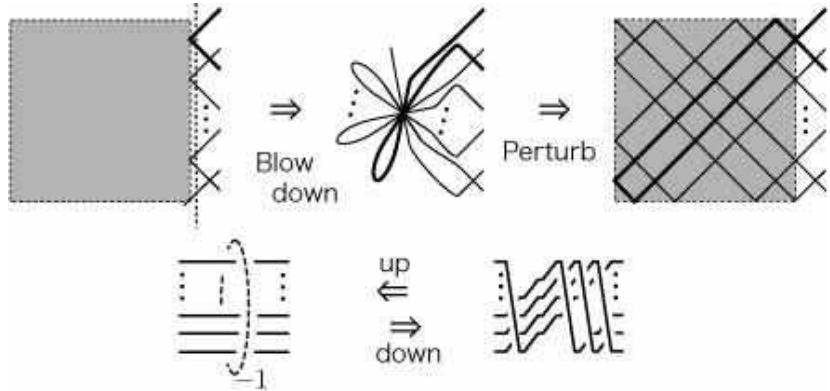}
\caption{Adding a square corresponds to a full-twist}
\label{fig:addS}
\end{center}
\end{figure}

Adding a square can be regarded as \lq\lq blow-down\rq\rq,
in the following sense.
The coordinate transformation $(x,y) = (x, xt)$ (or $=(yt, y)$) is called
a {\it blow up} in singularity theory, and is used for
resolution of singularity of complex curves, see \cite[p.16]{HKK}.
For example, for a coprime positive integers $(p,q)$ with $p < q$,
the complex curve $x^q - y^p =0$
becomes to $x^p(x^{q-p} - t^p) = 0$ by the transformation.
In this example, the link of the singularity at the origin changes
from the torus knot $T(p,q)$ to $T(p, q-p)$ and an unknot
defined by the complex line $t$-axis ($x=0$ in $xt$-plane) appears.
The unknot is the axis of the full-twist.
It corresponds to the Kirby calculus
in the bottom figure in Figure~\ref{fig:addS}, i.e.,
$(-1)$-framed unknot appears and the other components
change by a left-handed full-twist along the unknot.
Framings also change by a certain formula.
See \cite{K, GS} for Kirby calculus.

We call its inverse operation, i.e.,
the transformation $(x, y) = (x, y'/x)$ a blow down.
(For example, $y^2 = x + \epsilon$ becomes $y'^2 = x^2(x + \epsilon)$,
where both hand-sides are multiplied by $x^2$.)
Let $C$ be a real plane curve in $xy$-plane
that intersects with $y$-axis transversely.
By the transformation,
the all intersection points $C$ and $y$-axis concentrate to
the origin of $xy'$-plane,
and the left half ($x < 0$) of $C$ turns upside down
along $x$-axis, see the first arrow in Figure~\ref{fig:addS}.

For adding a square on an L-shaped curve along an edge of length $l$,
we first blow down the curve, where
we take $xy$-coordinate such that
$y$-axis is parallel and sufficiently close to the edge.
After that, we perturb the curve near the multiple crossing at the origin, 
see the second arrow in Figure~\ref{fig:addS}.
By using some $\Delta$-moves in Lemma~\ref{lem:divides}(6),
we can move the curve into
the required L-shaped curve of the square added L-shaped region.
The number of double points increases by $l(l-1)/2$.
For the Kirby calculus,
see the bottom figure in Figure~\ref{fig:addS} again.
By the operation, the knot $K$ changes by 
a right-handed full-twist along the unknot defined by the edge.
If $K$ is framed (i.e., with a surgery coefficient),
the framing increases by the square 
of the linking number of $K$ and $u$.
The linking number is equal to the length $l$ of the edge 
by Lemma~\ref{lem:divides}(3).

\subsection{Adding squares II}~\label{sbsec:AddS2}
\begin{figure}[h]
\begin{center}
\includegraphics[scale=0.4]{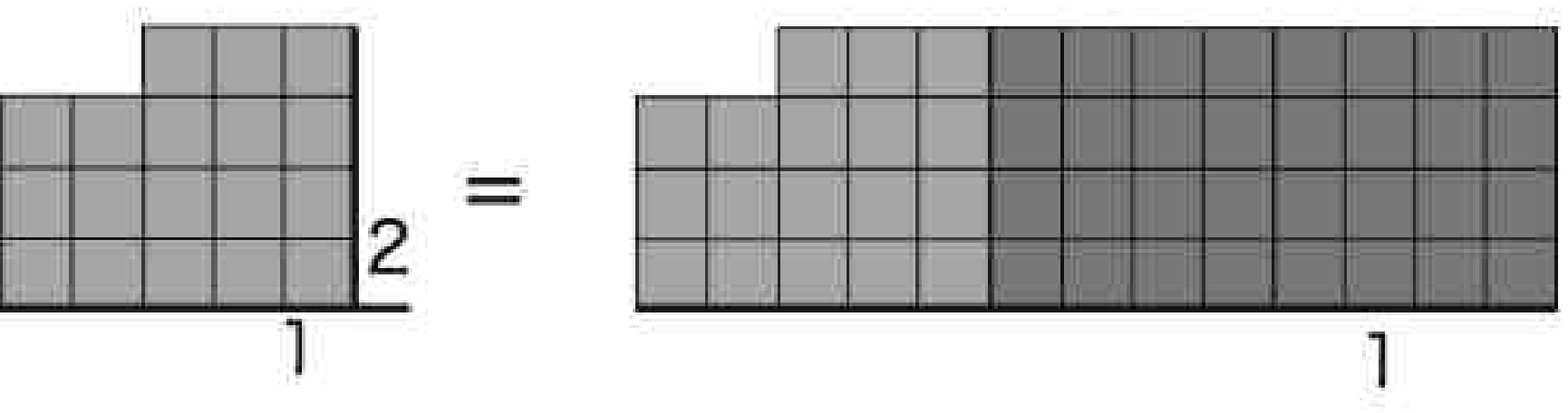}\\
\includegraphics[scale=0.4]{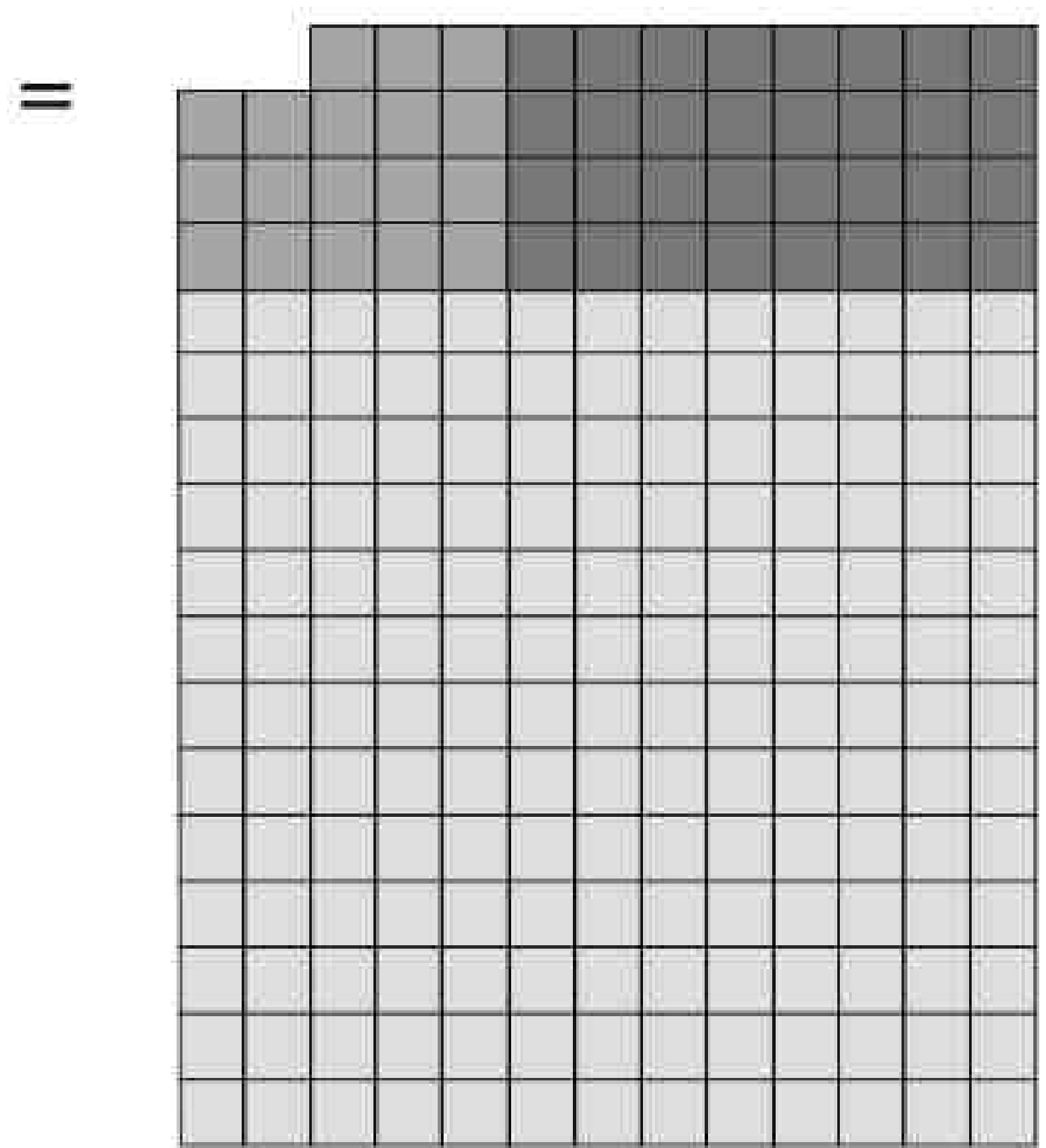}
\caption{Adding squares II}
\label{fig:LwithSq2}
\end{center}
\end{figure}

We can apply the operation adding squares twice
successively by changing the edges, as in Figure~\ref{fig:LwithSq2}.
In the same figure, we also define a drawing notation.
It is important which edge we apply the operation first.
\begin{figure}[h]
\begin{center}
\includegraphics[scale=0.5]{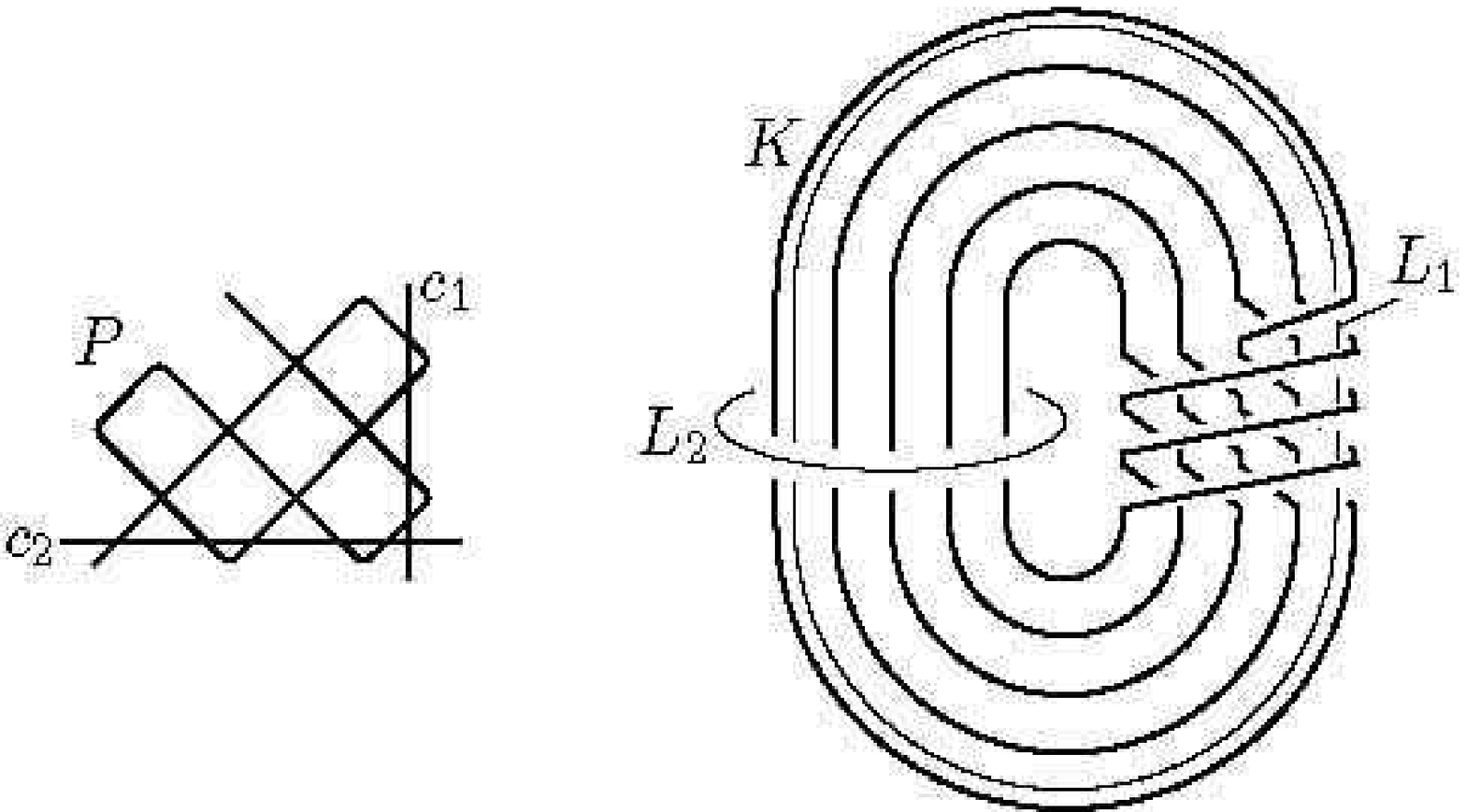}
\caption{The effect of adding squares twice}
\label{fig:Lcc}
\end{center}
\end{figure}

Here, we state the effect of twice adding squares
on the knots in $S^3$.
This is proved by Lemma~\ref{lem:pureAS}.
%
%
\begin{lem}~\label{lem:AddS2}
{\sl
Suppose that $P = X \cap \mathcal{L}$ be an L-shaped curve
and the first edge $e_1$ and the second one $e_2$ 
of the region $\mathcal{L}$ are specified.
By $P'' = X \cap \mathcal{L}''$,
we denote the resulting L-shaped curve obtained by
adding $n_1$ squares along $e_1$ first and
adding $n_2$ squares along $e_2$ second, successively.
Then, the divide knot $L(P'')$ is equal to the knot $K''$
obtained by two twistings from $L(P)$ in $S^3$ in the following sense:
\par
First, we take the three component divide link
$L(P \cup c_1 \cup c_2) = K \cup L_1 \cup L_2$
presented by the plane curve $P \cup c_1 \cup c_2$,
where $L(P)=K$ and $L_i$ is the component presented by slightly
pushed off $c_i$ of $e_i$ ($i=1,2$) into $\mathcal{L}$,
see Figure~\ref{fig:Lcc}.
Note that $L_1 \cup L_2$ is a Hopf link.
Next, we take $n_1$ full-twists of $K \cup L_2$ along $L_1$.
We call the resulting link $K' \cup L_2'$.
Finally, we take $n_2$ full-twists of $K'$ along $L_2'$.
We call the resulting knot $K''$.
}
\end{lem}
%
%
In Lemma~\ref{lem:AddS2}, $n_1, n_2$ are supposed to be positive,
however, regarding the statement as
{\it construction} of the knot $K''$ from $L(P)$ by two twistings,
it works also in the case $n_1, n_2 <0$.
The knot $K''$ may be no longer a divide knot,
by the obstruction of
braid (quasi-)positivity of divide knots in Lemma~\ref{lem:divides}(7)
(and Lemma~\ref{lem:QtwistT}).
In the next subsection, we consider the case $K''$ is
\underline{the mirror image} of an L-shaped divide knot.

\subsection{Adding squares III}~\label{sbsec:AddS3}
We extend the operation adding $n$ squares
into the case $n <0$ partially,
in analogy with Lemma~\ref{lem:AddS2}.
It corresponds to negative (i.e., left-handed) $\vert n\vert$ full-twists,
and will be called \lq\lq adding negative squares\rq\rq.
We consider the case that
the resulting knot of negative full-twists of an L-shaped divide knot
is the mirror image of another L-shaped divide knot
under some conditions.
%
%
\begin{defn}~\label{def:AddNS}
{\rm (Adding negative squares in a certain case) \ }
{\rm
Let ${\mathcal L}$ be
an L-shaped region of type $[a_1, a_2; b_1, b_2]$
with a specified edge $e$.
We assume that the edge is the (bottom) one of length $a_2$,
\[
b_2 = b_1+1
\quad \textrm{and} \quad
\vert n \vert a_2 > b_1+1.
\] 
Only under this condition,
we define {\it adding $n$ squares with $n<0$ along $e$}
as that the resulting region is
of type
\[
[a_2-a_1+1, a_2; \vert n \vert a_2-b_1-1, \vert n \vert a_2-b_1 ].
\]
}
\end{defn}
%
%
\begin{figure}[h]
\begin{center}
\includegraphics[scale=0.5]{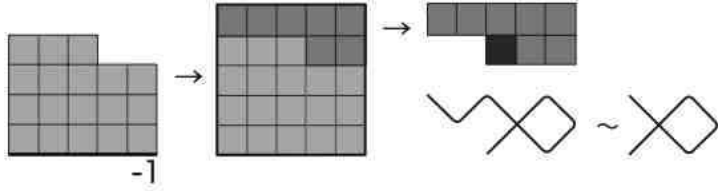}
\caption{Adding negative squares}
\label{fig:LwithSq3}
\end{center}
\end{figure}
This operation should be regarded geometrically as follows,
see Figure~\ref{fig:LwithSq3}:
We assume that the initial L-shaped region is at the origin as
in Definition~\ref{def:L-sh} (once forgetting the assumption ($\ast$))
to explain the operation by using $xy$-coordinate.
Then adding $n$ squares with $n<0$ is defined as
\begin{eqnarray*}
& & 
\textrm{cl}
\big [
\{
(x,y)
\vert \
0 \leq x \leq a_2, \
0 \leq y \leq \vert n \vert a_2
\}
\setminus
{\mathcal L}
\big ]
\\
& & \qquad \cup \
\{
(x,y)
\vert \
a_1 -1 \leq x \leq a_1, \
b_1 \leq y \leq b_1 +1
\}.
\end{eqnarray*}
By the finally added unit square at the concave point,
if the concave point of the initial region is at an odd point,
then that of the resulting region is also at an odd point,
i.e., we can keep the condition ($\ast$).
%
%
\begin{lem}~\label{lem:AddNS}
{\sl
Under the condition of adding negative squares in
Definition~\ref{def:AddNS},
adding $n$ squares with $n<0$ on an
L-shaped curve $P$ along the edge
corresponds to taking the mirror image of
$n$ right-handed (i.e., $\vert n \vert$ left-handed) full-twists
on the divide knot $L(P)$ along the unknot defined by the edge.
}
\end{lem}
%
%
{\it Proof.}
\quad 
By the operation,
the type of L-shaped curves is changed
from $[a_1, a_2; b_1, b_1 +1]$ to
$[a_2-a_1+1, a_2; \vert n \vert a_2-b_1-1, \vert n \vert a_2-b_1]$.
By Lemma~\ref{lem:QtwistT}, the initial curve
presents the closure of the braid $W(a_2)^{b_1}W(a_1)$
of index $a_2$.
The edge presents the braid axis in $S^3$,
The positive full-twist is $W(a_2)^{a_2}$ in this situation,
and is in the center of the braid group.
The resulting knot of the $n$ full-twists is
the closure of
$W(a_2)^{b_1+ na_2}W(a_1)
= W(a_2)^{-(\vert n \vert a_2-b_1)}W(a_1)$.
Note that $\vert n \vert a_2 > b_1+1$ is assumed.
Lemma~\ref{lem:AddNS} follows from
Lemma~\ref{lem:mirror}($-+$).
\qed
%
%
\begin{ex}~\label{ex:237b}
{\rm
The example (from $[3,5;3,4]$ to $[3,5;1,2]$) 
in Figure~\ref{fig:LwithSq3}
shows that the mirror image of $T(2,3)$
is obtained by $P(-2,3,7)$ (in Example~\ref{ex:237})
by ($-1$) full-twist along the unknot defined by the bottom edge,
whose linking number with $P(-2,3,7)$ is $\pm 5$.
Note that L-shaped curve of type $[3,5;1,2]$ is the same curve with 
the billiard curve $B(2,3)$ defined in Lemma~\ref{lem:GHY}.
}
\end{ex}
%
%
\begin{quest}
{\rm
Extend the operation adding negative squares into (more) general cases.
}\end{quest}
%
%
\par

\subsection{How to construct the L-shaped curve}~\label{sbsec:demo}
Preparation is complete.
For each Berge's knot
$K = K_{{\mathcal X}}(\delta, \varepsilon, A, k, t)$
in Type ${\mathcal X}$, we take the L-shaped region
${\mathcal L} = {\mathcal L}(\varepsilon, A, k, t)$
in Table~\ref{tbl:Lregion},
where we used the drawing notation of adding squares.
Then, the plane curve $P = X \cap {\mathcal L}$, as a divide,
presents the knot $K$ or its mirror image.
In fact, each L-shaped region in Table~\ref{tbl:Lregion}
is carefully constructed such that
Berge's braid presentation of the knot in Lemma~\ref{lem:Berge}
agrees with the braid presentation of the region in Lemma~\ref{lem:mirror}
under the suitable choice of $\delta$.
The proof is in the next section. 
\par
In Table~\ref{tbl:Lregion}, 
we draw each L-shaped region in the case of the smallest $A$.
\par
\medskip
\noindent
{\bf Demonstration.}
$K_{\textrm{III}}(\delta, \varepsilon, A, k, t)$
for $(\delta, \varepsilon, A, k, t) = (1,1, 2, 2,1)$. 
\par
By Table~\ref{tbl:Para}(1),
$B = A(3+2k)-\varepsilon = 13, \
b = - \delta \varepsilon (2A + tB )= -17$.
By Lemma~\ref{lem:Berge},
it has a braid presentation $W(13)^{-17} W(3)^{1}$.
The surgery coefficient is $bB + \delta A = -219$.
On the other hand, according to Table~\ref{tbl:Lregion},
the L-shaped region $\mathcal{L}(\varepsilon, A, k, t)$
with $(\varepsilon, A, k, t) = (1, 2, 2,1)$ is
the L-shaped region of type $[11,13; 16,17]$
(the region at the bottom in Figure~\ref{fig:LwithSq2}),
whose area is $219$.
By Lemma~\ref{lem:QtwistT}, its corresponding plane curve
presents the closure of $W(13)^{16} W(11)$.
By Lemma~\ref{lem:Wb-1}, the knot is equal to the closure
of $W(13)^{17} W(3)^{-1}$, which is the mirror image
of the required knot.
\par

\section{Proof of Theorem~\ref{thm:main1} and \ref{thm:main2}}~\label{sec:Proof}
Theorem~\ref{thm:main1} is proved by verifying that
Berge's braid presentation of the knot (Lemma~\ref{lem:Berge})
and that of the L-shaped region in Table~\ref{tbl:Lregion}
(Lemma~\ref{lem:mirror}) agree,
under the suitable choice of the sign $\delta$ in each Type.
In the proof below, we will also decide 
the choice of $\delta$ (depending on $\mathcal{X}$ and $(\varepsilon, A, k, t)$).
We denote the result by
$\delta_{\mathcal{X}}(\varepsilon, A, k, t)$, or $\delta_{\mathcal{X}}$ for short.
It will be used in the proof of Lemma~\ref{lem:pos}
\par \medskip 
\noindent
{\it Proof.} (of Theorem~\ref{thm:main1})
\quad 
Here, we prove Theorem~\ref{thm:main1} only in the case of Type III.
The proofs in the other cases are similar.

In Type III, in Table~\ref{tbl:Para}(1) and (2), we find
\begin{eqnarray*}
a=0, \quad
B = Al-\varepsilon, \quad
l = 3 +2k, \quad
b = - \delta \varepsilon (2A + t B).
\end{eqnarray*}
Thus, Berge's braid presentation $W(B)^bW(A+1-a)^{\delta}$ of the knot
$K_{\textrm{III}}(\delta, \varepsilon, A, k, t)$ in Lemma~\ref{lem:Berge}
(\cite{Bg2}) is
%
\begin{eqnarray}~\label{eq:IIIa}
W((3+2k)A-\varepsilon)^{- \delta \varepsilon (2A + t B)}
\, W(A+1)^\delta
\end{eqnarray}
%
First, we consider the case $k=t=0$.
\par
\noindent
{\bf Case 1+} \ (Type III, $\varepsilon =1, k=t=0$) \\
In Table~\ref{tbl:Lregion}, we find that the L-shaped region is of type
$[2A-1, 3A-1; 2A-1, 2A]$.
On the other hand, Berge's presentation (\ref{eq:IIIa}) is now
$W(3A-1)^{-2\delta A} W(A+1)^{\delta}$.
Here we choose $\delta = -1$.
We use Lemma~\ref{lem:mirror}($+-$) on
the knot of $W(3A-1)^{2A} W(A+1)^{-1}$.
\par
\noindent
{\bf Case 1-} \ (Type III, $\varepsilon = -1, k=t=0$) \\
In Table~\ref{tbl:Lregion}, we find that the L-shaped region is of type
$[A+1, 3A+1; 2A, 2A+1]$.
On the other hand, Berge's presentation (\ref{eq:IIIa}) is now
$W(3A+1)^{2\delta A} W(A+1)^{\delta}$.
We choose $\delta = 1$. We use Lemma~\ref{lem:mirror}($++$).
\par \medskip
Next, we consider the case $k > 0$ ($t= 0$).
In this case, we use  the symmetry of Corollary~\ref{cor:sym}
to verify that the parameter $k$ contributes as the $k$ full-twists
on the knots.
\par \medskip
\noindent
{\bf Case 2+} \ (Type III, $\varepsilon = 1, k >0, t=0$) \\
The L-shaped region in Table~\ref{tbl:Lregion} is of type
$[2A -1 + 2kA, 3A -1 + 2kA; 2A-1, 2A]$,
which presents the closure of
$W(3A -1 + 2kA)^{2A-1} W(2A -1 + 2kA)$ by Lemma~\ref{lem:QtwistT}.
By the symmetry of Corollary~\ref{cor:sym}, it presents
the same knot of
\[
W(2A)^{2A-1 + 2kA} W(2A -1)^A
=
(W(2A)^{2A})^k \, W(2A)^{2A-1} W(2A-1)^A.
\]
Its first part $(W(2A)^{2A})^k$ means $k$ full-twists.
On the other hand, Berge's presentation (\ref{eq:IIIa}) is now
$W(3A -1 + 2kA)^{-2\delta A} W(A+1)^{\delta}$
We choose $\delta =-1$.
We use Lemma~\ref{lem:mirror}($+-$).
\par
\noindent
{\bf Case 2-} \ (Type III, $\varepsilon = -1, k >0, t=0$) \\
The L-shaped region in Table~\ref{tbl:Lregion} is of type
$[A+1, 3A+1+ 2kA; 2A, 2A+1]$, which presents the closure of
$W(3A+1+ 2kA)^{2A}W(A+1)$ by Lemma~\ref{lem:QtwistT}.
By Corollary~\ref{cor:sym}, it presents the same knot of
\[
W(2A+1)^{A+1} W(2A)^{A + 2kA}
=
W(2A+1)^{A+1} W(2A)^A \, (W(2A)^{2A})^k.
\]
Its final part $(W(2A)^{2A})^k$ means $k$ full-twists
of $2A$ strings of the braid of index $2A+1$.
On the other hand, Berge's presentation (\ref{eq:IIIa}) is now
$W(3A +1 + 2kA)^{2\delta A} W(A+1)^{\delta}$ by (\ref{eq:IIIa}).
We choose $\delta =1$. We use Lemma~\ref{lem:mirror}($++$).
\par \medskip
Before we go into the case $t \not= 0$, we remark the followings.
\begin{enumerate}
\item[(1)] The parameter $t$ can be negative.
\item[(2)] The parameters $A, k$ and $B$ are
independent from $t$.
\item[(3)] 
The parameter $B$
is equal to the length of the edge that is added $t$ squares.
\item[(4)] In every case in Table~\ref{tbl:Lregion},
we can apply the operation adding $t$ squares along the edge,
i.e., the condition
\lq\lq $b_2 = b_1 +1$ and $\vert t \vert a_2 > b_1+1$\rq\rq \ %
for adding negative squares are satisfied, even if $t=-1$.
\end{enumerate}
From now on, we use
Berge's braid presentation of 
$K_{\textrm{III}}(\delta, \varepsilon, A, k, t)$
in Lemma~\ref{lem:Berge} (\cite{Bg2}) in the form
%
\begin{eqnarray}~\label{eq:IIIb}
W(B)^{-\delta \varepsilon (2A+tB)} W(A+1)^{\delta}.
\end{eqnarray}
%
\par
\noindent
{\bf Case 3+} \ (Type III, $\varepsilon = 1, t > 0$) \\
The L-shaped region in Table~\ref{tbl:Lregion} is of type
$[B-A, B; 2A -1+ tB, 2A+tB] $ with $B = 3A -1+2kA$.
On the other hand, Berge's presentation (\ref{eq:IIIb}) is now
$W(B)^{-\delta (2A+tB)} W(A+1)^{\delta}$
We choose $\delta = -1$.
We use Lemma~\ref{lem:mirror}($+-$).
\par
\noindent
{\bf Case 3-} \ (Type III, $\varepsilon = -1, t > 0$) \\
The L-shaped region in Table~\ref{tbl:Lregion} is of type
$[A+1, B; 2A + tB, 2A +1+ tB ] $ with $B = 3A +1+2kA$.
On the other hand, Berge's presentation (\ref{eq:IIIb}) is now
$W(B)^{\delta (2A+tB)} W(A+1)^{\delta}$.
We choose $\delta =1$. We use Lemma~\ref{lem:mirror}($++$).
\par \medskip
Before we consider the case $t<0$, we remark
the condition $2A \leq B$ in Section~\ref{sec:Berge}.
In fact, $2A = B$ never occurs.
\par \medskip
\noindent
{\bf Case 4+} \ (Type III, $\varepsilon = 1, t < 0$) \
We use adding negative squares. \\
The type of L-shaped region in Table~\ref{tbl:Lregion} is
$[A+1, B; \vert t\vert B-2A, \vert t\vert B-2A + 1]$
with $B = 3A -1+2kA$.
On the other hand, Berge's presentation (\ref{eq:IIIb}) is now
$W(B)^{-\delta (2A+tB)} W(A+1)^{\delta}$.
We choose $\delta =1$. We use Lemma~\ref{lem:mirror}($++$)
on $W(B)^{\vert t \vert B -2A} W(A+1)$.
\par
\noindent
{\bf Case 4-} \ (Type III, $\varepsilon = -1, t < 0$) \
We use adding negative squares. \\
The type of L-shaped region in Table~\ref{tbl:Lregion} is
$[B-A, B; \vert t\vert B-2A -1, \vert t\vert B-2A] $ with $B = 3A +1+2kA$,
which presents the closure of
$W(B)^{\vert t\vert B-2A-1} W(B-A)$.
On the other hand, Berge's presentation (\ref{eq:IIIb}) is now
$W(B)^{\delta (2A+tB)} W(A+1)^{\delta}$
We choose $\delta = -1$.
We use Lemma~\ref{lem:mirror}($+-$) on
the knot of $W(B)^{\vert t \vert B -2A}  W(A+1)^{-1}$.
\par
After all, for given parameters $(\varepsilon, A, k, t)$,
we have shown that
the L-shaped curve $P = X \cap \mathcal{L}$ of the
L-shaped region
$\mathcal{L} = \mathcal{L}_{\textrm{III}}(\varepsilon, A, k, t)$
in Table~\ref{tbl:Lregion}
presents $K_{\textrm{III}}(\delta, \varepsilon, A, k, t)$
for suitable choice of $\delta$.

The proof of Theorem~\ref{thm:main1} in the case of Type III is completed.
The cases of the other Types are proved by the same argument.
\qed
\par \medskip

In the proof of Theorem~\ref{thm:main1},
$\delta_{\textrm{III}}(\varepsilon, A, k, t)$,
the suitable choice of $\delta$ in Type III
is determined. Considering the other Types, it extends as
\[
\delta_{\mathcal{X}}(\varepsilon, A, k, t) =
\begin{cases}
-1 & \textrm{ if } \varepsilon \cdot \text{sgn}(t)= +1 \\
 1 & \textrm{ if } \varepsilon \cdot \text{sgn}(t)= -1 \\
\end{cases},
\]
where $\text{sgn}(t) = 1$ if $t >0$ or $t=0$, $\text{sgn}(t) = -1$ otherwise.
It means simply
%
\begin{eqnarray}~\label{eq:det}
-\delta_{\mathcal{X}} \cdot \varepsilon \cdot \text{sgn}(t)= +1.
\end{eqnarray}
%
Using $\delta_{\mathcal{X}}$, Theorem~\ref{thm:main1} means simply
\[
L(\, P_{\mathcal{X}}(\varepsilon, A, k, t) \, )
=
K_{\mathcal{X}}(\delta_{\mathcal{X}}, \varepsilon, A, k, t).
\]

Now, we consider the coefficients of the lens space surgeries.
Before we start the proof of Theorem~\ref{thm:main2}, we show
%
%
\begin{lem}~\label{lem:pos}
{\sl
The surgery coefficient of
$L(\, P_{\mathcal{X}}(\varepsilon, A, k, t) \, )$
is positive.
}
\end{lem}
%
%
{\it Proof.}
\quad 
As in the proof of Theorem~\ref{thm:main1}, here
we prove the lemma only in the case of Type III.
In Table~\ref{tbl:Para}(3), we find that the coefficient of
$L(\, P_{\textrm{III}}(\varepsilon, A, k, t) \, ) =
K_{\textrm{III}}(\delta_{\mathcal{X}}, \varepsilon, A, k, t)$
(as in Type III) is
\[
-\delta_{\mathcal{X}} \cdot \varepsilon \cdot
(6A^2- 3\varepsilon A + k(2A)^2 + t B^2).
\]
Since $\vert B \vert = (3+2k)A - \varepsilon$
(the length of the longest edge), it holds that
$\text{sgn}(6A^2- 3\varepsilon A + k(2A)^2 + t B^2)  = \text{sgn}(t)$.
Thus the sign of the coefficient is
equal to
$-\delta_{\mathcal{X}}\cdot \varepsilon \cdot \text{sgn}(t) = +1$
by (\ref{eq:det}).
The proof of the other Types are similar.
\qed
%
%
\begin{quest}
{\rm
Prove Lemma~\ref{lem:pos} without Berge's braid presentation.
If a knot $K$ is a closure of a positive braid and 
admits a lens space surgery, is the coefficient positive?
}\end{quest}
%
%

Now, we prove a precise version of Theorem~\ref{thm:main2} on
the difference between the area
$\textrm{area}(P)$ of the L-shaped curve
and the surgery coefficient
$\textrm{coef}(L(P))$ ($>0$)
of the lens space surgery as in Type $\mathcal{X}$.
%
%
\begin{lem}~\label{lem:v-c}
{\sl 
Under the correspondence in Theorem~\ref{thm:main1},
it holds that
\[
\textrm{area} (P) - \textrm{coef}(L(P)) =
\begin{cases}
0 & \textrm{ if } (-1)^a \cdot \varepsilon \cdot \text{sgn}(t)= +1 \\
1 & \textrm{ if } (-1)^a \cdot \varepsilon \cdot \text{sgn}(t)= -1 \\
\end{cases}.
\]
where $\text{sgn}(t) = 1$ if $t >0$ or $t=0$, $\text{sgn}(t) = -1$ otherwise.
See Table~\ref{tbl:Para}(1) for the definition of $a$ ($=0$ or $1$).
}
\end{lem}
%
%
{\it Proof.}
\quad 
First, in the case $k = t = 0$, it is easy to verify the equation
in Table~\ref{tbl:Lregion}.
The parameters $k$ and $t$ with $t>0$ contribute as
adding positive squares.
In the operation adding a positive square along an edge of length $x$,
the area increases by $x^2$. 
It is compatible with the terms $+ kA^2 + tB^2$ (or $+ k(2A)^2 + tB^2$)
in the surgery coefficients in Table~\ref{tbl:Para}(3).
In the case $t < 0$, we do the operation adding negative squares.
Suppose that we get the curve $P_{\textrm{new}}$ from $P_{\textrm{old}}$
by adding $t$ squares with $t<0$ along an edge of length $x$.
Then
\[
\textrm{coef}(L(P_{\textrm{new}}))
= 
- ( \, \textrm{coef}(L(P_{\textrm{old}})) - \vert t \vert x^2 \,),
\]
since the new divide knot $L(P_{\textrm{new}})$ is the mirror image of the
knot obtained by left-handed $\vert t \vert$ twists from $L(P_{\textrm{old}})$.
On the other hand,
\[
\textrm{area}(P_{\textrm{new}})
=
\vert t \vert x^2 - \textrm{area}(P_{\textrm{old}}) +1,
\]
where the last $+1$ corresponds to the finally added unit square.
Thus
\[
\textrm{area} (P_{\textrm{new}}) - \textrm{coef}(L(P_{\textrm{new}}))
=
1 - ( \,
\textrm{area} (P_{\textrm{old}}) - \textrm{coef}(L(P_{\textrm{old}})) \, ).
\]
We have the lemma.
\qed
\par \medskip

Note that Berge's constant $a$ in \cite[Table 3(p.15)]{Bg2} was defined
geometrically in the context of doubly-primitive knots.
It might be curious that $a$ is related to the
difference between the area and the coefficient
as above.
\par \medskip

Goda and Teragaito \cite{GT} conjectured
an inequality
\[
2g(K)+ 8 \leq \vert r \vert  \leq 4g(K) -1,
\]
on the surgery coefficient $r$
and the genus $g(K)$ of the hyperbolic lens space surgery $(K,r)$.
It is called \lq\lq Goda--Teragaito conjecture\rq\rq.
We are concerned with the left hand-side inequality.
By Proposition~\ref{prop:GT} and Lemma~\ref{lem:v-c},
we have:
%
%
\begin{cor}
{\sl 
Let $[a_1, a_2; b_1, b_2]$ be
the type of the L-shaped curve of
our presentation (in Table~\ref{tbl:Lregion})
of Berge's knot in Type III, IV, V and VI as a divide knot $L(P)$.
Then, it holds that
\[
\textrm{coef}(L(P)) - 2g(L(P)) = a_2 +b_2 -1, \
\textrm{or } \ a_2 + b_2 -2.
\]
}
\end{cor}
%
%
The parameters $A, k$ or $\vert t\vert$ are greater,
the difference $\vert r \vert - 2g(K)$ can be greater.
\par

\section{Further Observation}~\label{sec:Obs}
Divide presentation of L-shaped divide knots
helps us to study the constructions of the knots, and
the relationship among the lens space surgeries.

\subsection{Twisted torus knots}~\label{sbsec:TwT}
Following Dean \cite{De},
by {\it a twisted torus knot} $T(p,q; r,s)$, if $r<p$,
we mean the knot obtained from the torus knot
$T(p,q)$ by $s$ full-twists of $r$ strings
in the $p$ parallel strings of $T(p,q)$ in the standard position.
On the other hand, if $r>p$, we mean the knot obtained from
the torus knot $T(p,q)$ as a closure of the braid
$w(p,q) \sigma_p \sigma_{p+1} \cdots \sigma_{r-1}$
by $s$ full-twists of all $r$ strings,
where $w(p,q) \sigma_p \sigma_{p+1} \cdots \sigma_{r-1}$ is
a positive Markov stabilization of
the standard braid $w(p,q)$ of $T(p,q)$ of index $p$
to a braid of index $r$.
The following lemma follows from the braid presentation
in Lemma~\ref{lem:QtwistT} (and Lemma~\ref{lem:pureAS}).
\par
\begin{figure}[h]
\begin{center}
\includegraphics[scale=0.4]{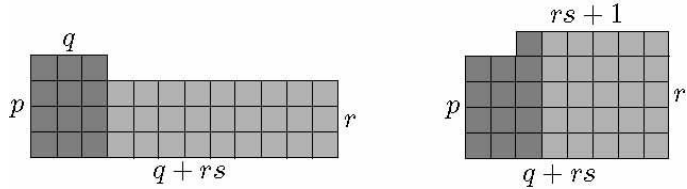}
\caption{Twisted torus knots $T(p,q;r,s)$ \ (ex. $T(4, 3 ; 3, 3)$ and $T(4,3; 5,1)$)}
\label{fig:TwTL}
\end{center}
\end{figure}
\begin{figure}[h]
\begin{center}
\includegraphics[scale=0.5]{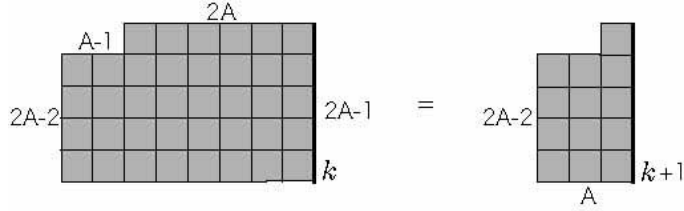}
\caption{$K_{\textrm{V}}(-1,-1, A, k, -1) =  T(2A-2, A; 2A-1, k+1)$
}
\label{fig:TwT}
\end{center}
\end{figure}

%
%
\begin{lem}~\label{lem:TwTL}
{\sl
Let $(p,q)$ be a coprime pair of positive integers, and
$r, s$ integers satisfying $0 < r \not= p$ and $s > 0$.
The twisted torus knot $T(p,q; r,s)$
is one of A'Campo's divide knots, and can be presented
by an L-shaped curve of type
\[
\begin{cases} \
[q, q+rs; r, p]  & \textrm{ if $r < p$} \\
[rs+1, q+rs; p, r ]  & \textrm{ if $r > p$}
\end{cases},
\]
as A'Campo's divide knots,
see Figure~\ref{fig:TwTL}.
}
\end{lem}
%
%
Note that, a twisted torus knot $T(p,q; r,s)$ can be
accidentally non-hyperbolic,
a torus knot $T(p',q')$, or a cable knot $C( T(p',q'); m', n')$
of a torus knot, see \cite{Y2, MY} for such phenomena.
%
%
\begin{lem}~\label{lem:TwT}
{\sl
Each knot in the following list is a twisted torus knot: 
\par
(knots with $t=0$) \par
\begin{tabular}{ll}
$K_{\textrm{III}}(1,-1, A, k, 0)$
&
$= T(2A+1,A+1; 2A, k+1)$,
\\
$K_{\textrm{IV}}(-1,1, A, k, 0)$
&
$= T(A, kA + (3A+1)/2; A-1, 1)$,
\\
$K_{\textrm{V}}(-1,1, A, k, 0)$
&
$= T(A, (k+1)A+2; A-1, 1)$,
\\
$K_{\textrm{VI}}(1,-1, A, 0, 0)$
&
$= T(A-1, A+1; A, 1)$,
\\
\end{tabular} \par
(knots with $t=-1$) \par
\begin{tabular}{ll}
$K_{\textrm{III}}(1,1, A, k, -1)$
&
$= T(A, A +1 ; A-1, k+2 )$,
\\
$K_{\textrm{IV}}(-1,-1, A, k, -1)$
&
$= T((3A-1)/2, A; (3A+1)/2, k+1)$,
\\
$K_{\textrm{V}}(-1,-1, A, k, -1)$
&
$= T(2A-2, A; 2A-1, k+1)$.
\\
\end{tabular}
}
\end{lem}
%
%
{\it Proof.}
\quad 
Only we have to do is to verify them by comparing
types of the L-shaped regions,
thus we omit the proof in detail.
See Figure~\ref{fig:TwT} for the case of
$K_{\textrm{V}}(-1,-1, A, k, -1)$, which needs
adding $t$ square with $t = -1$.
\qed \par
\medskip

Now, we count $s$ full-twists along a fixed unknot as one twisting.
If after $s_1$ full-twists along an unknot we take another
$s_2$ full-twists along a different unknot, then we count the operation as
two twistings.
%
%
\begin{cor}~\label{cor:TwT}
{\rm ([DMM])} \quad 
{\sl
Every knot in Berge's in Type III, IV, V and VI is obtained
at most one ($s'$) full-twists from a twisted torus knot 
$T_w := T(p,q; r,s)$, thus 
is obtained at most two twistings from a torus knot, such that 
every knot in the twisting process 
($T(p,q; r,i)$ with $0 \leq i \leq s$ and 
$j$ full-twists of $T_w$ with $0 \leq j \leq s'$) admits lens space surgery.
Furthermore, as a twisted torus knot above, we can take 
$T(p,q; r,s)$ that satisfies $\vert r - p \vert =1$. 
}
\end{cor}
%
%
In the process of preparation of this paper, the author was
informed by Motegi (\cite{DMM}) that every knot in Berge's list
(including Type VII, ..., XII, see Section~\ref{sec:intro}) is 
obtained by at most two twistings from a torus knot,
which includes the above Corollary. 
\par

\subsection{Relations between different Types}~\label{sbsec:type}
\begin{figure}[h]
\begin{center}
\includegraphics[scale=0.5]{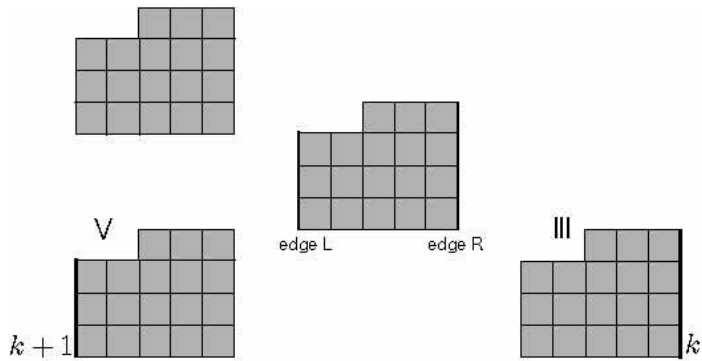}
\caption{Type III and Type V}
\label{fig:Type35}
\end{center}
\end{figure}

See the L-shaped region in the left top figure in Figure~\ref{fig:Type35}.
It is of type $[3,5; 3,4]$, and presents $K_{\textrm{III}}(-1,1, 2, 0, 0)$
($=P(-2,3,7)$).
The L-shaped regions obtained by adding $k$ squares
along the right edge (denoted by R),
we have a subsequence
$K_{\textrm{III}}(-1,1, 2, k, 0)$ in Type III,
see Type III($\varepsilon =1$) in Table~\ref{tbl:Lregion}.
On the other hand, those obtained by adding $k+1$ squares
with $k \geq 0$ along the left edge (denoted by L),
we have another subsequence
$K_{\textrm{V}}(1,1, 3, k, 0)$ in Type V (a different Type),
see Type V($\varepsilon = -1$) in Table~\ref{tbl:Lregion}.
By
\[
K_{\textrm{III}}(-1,1, 2, 0, 0) \mapsto
K_{\textrm{V}}(1,1, 3, k, 0),
\]
we denote such a relation, regarding it as
$K_{\textrm{III}}(-1,1, 2, 0, 0)=$
\lq\lq$K_{\textrm{V}}(1,1, 3, -1, 0)$\rq\rq.
Using L-shaped curve presentation in Table~\ref{tbl:Lregion},
we can see such relations more: 
\par
$\cdot$ For each $A$ with $A \geq 3$,
$K_{\textrm{III}}(-1,1, A, 0, 0)
\mapsto
K_{\textrm{IV}}(1,-1, 2A-1, k, 0)$. 
\par
$\cdot$ For each $A$,
$K_{\textrm{III}}(1,-1, A, 0, 0)
\mapsto
K_{\textrm{IV}}(-1,1, 2A+1, k, 0)$.
\par
%
%
\begin{quest}
{\rm
Find such relations more, especially in the case $t\not=0$.
}\end{quest}
%
%

\subsection{Parameters Translation}~\label{sbsec:parameter}
Parameters of Berge's knots are different among some papers.
In Table~\ref{tbl:Para}(4),
we give the translation formula between
\begin{center}
$(A, k)$ in \cite{Bg2} and this paper
\quad and \quad
$(n, p)$ in \cite{Bg, Ba, Ba3, DMM}.
\end{center}
The signs $\delta, \varepsilon$ are commonly used
among these papers.

For example, $K_{\textrm{IV}}(\delta ,-1, 5, 1, t)$
(i.e., $( \varepsilon, A, k) = (-1, 5, 1)$) in this paper
is, up to mirror image, obtained
from the knot of $(\varepsilon, n, p) = (-1, 2, 4)$ in Type IV
in \cite{Bg, Ba, Ba3, DMM} by $\pm t$ full-twists.

The parameter $p$ is defined as a positive integer, and
our $k$ is just a parallel shift of $p$ such that
\lq\lq $k=0$ at the minimal possible value as $p$\rq\rq:
the statement
\lq\lq $\varepsilon p \not= -2, -1, 0, 1$\rq\rq \ %
(see Table~\ref{tbl:Para}(4)) is referred as
\begin{center}
\lq\lq $p \geq 2$ if $\varepsilon = +1$ and
$p \geq 3$ if $\varepsilon = -1$\rq\rq \ %
\quad 
in \cite{Bg, Ba, Ba3, DMM}.
\end{center}
On the other hand, we say
\begin{center}
\lq\lq $k := p-2$ if $\varepsilon = +1$ and
$k := p-3$ if $\varepsilon = -1$\rq\rq,
\quad 
in the present paper.
\end{center}
\par

\bigskip
{\bf Acknowledgement.}
The author would like to thank to 
Professor Mikami Hirasawa, 
Dr. Tomomi Kawamura,
Dr. Masaharu Ishikawa,
Professor Sergei Chmutov,
and Professor Norbert A'Campo
for informing him on A'Campo's divide knot thoery.
The author also would like to thank to 
Professor Kimihiko Motegi, 
Professor Masakazu Teragaito, 
Professor Hiroshi Goda, 
Professor Noriko Maruyama,
Dr. Toshio Saito,
Dr. Kenneth Baker, 
Dr. Arnaud Deruelle,
Dr. Hiroshi Matsuda,
and 
Professor John Berge 
for helpful suggestion on lens space surgery.
\bigskip



{\small
\par
Dept. of Systems Engineering, The Univ. of Electro-Communications \par
1-5-1,Chofugaoka, Chofu, Tokyo, 182-8585, JAPAN \par
}
{\tt yyyamada@sugaku.e-one.uec.ac.jp} \par

%
%


\bigskip

%
%
\begin{table}[h]
\begin{center}
\begin{tabular}{l|l|l|l|lll}
&
$a$ & $A$ & $l$ & & 
$\delta, \varepsilon \in \{ \pm 1 \}$, \ $t \in {\bf Z}$ \\
\hline
{\bf III} 
& $0$ &
$2, 3, 4, \cdots$
& 
$3, 5, 7, \cdots$(odd)
&
$B = Al -\varepsilon$,
&
$b = - \delta \varepsilon (2A + tB)$ \\
{\bf IV}
& $1$ &
$5, 7, 9, \cdots$(odd)
& 
$5, 7, 9, \cdots$(odd)
&
$B = (Al -\varepsilon)/2$,
&
$b = - \delta \varepsilon (A + tB)$ \\
{\bf V}
& $1$ &
$3, 5, 7, \cdots$(odd)
& 
$2, 3, 4, \cdots$
&
$B = Al + \varepsilon^*$,
&
$b = - \delta \varepsilon (A + tB)$ \\
{\bf VI}
& $0$ &
$4, 6, 8, \cdots$(even)
& 
&
$B = 2A+1$,
&
$b = \delta (A-1 + tB)$ \\
\hline
\end{tabular} \\
\hfill {\small Here $\varepsilon^* \not = -1$ if $l=2$ (since $0 <2A \leq B$).}
\par \medskip
(1) The parities and ranges of $A, l$ and the settings of $a, B, b$ (\cite{Bg2})
\par \vskip20pt
%
%
\begin{tabular}{lcrcl}
{\bf III}: \  $l = 3+ 2k $, & & 
{\bf V}: \ $l = 2+ k$ \ (if $\varepsilon =+1$),
& & 
{\bf VI}: \ $k \equiv 0$ \ ($\varepsilon \equiv -1$). \\
{\bf IV}: \  $l = 5+ 2k $, & & 
 $l = 3+ k$ \ (if $\varepsilon = -1$), 
\end{tabular}
\par \medskip
(2) Parameter $k$
\par \vskip20pt
\begin{tabular}{ll}
{\bf III.}
& $
- \delta \varepsilon 
\left (6 A^2 - 3 \varepsilon A + k(2A)^2 + t B^2 \right )$, \\
{\bf IV.}
& $ - \delta \varepsilon \left ( \dfrac52 A^2 - \dfrac32 \varepsilon  A
+ k A^2 + t B^2 \right )$, \\
{\bf V.}
& $ - \delta  \left ( 2A^2 + k A^2 + t B^2 \right )$ \ if $\varepsilon =+1$,
\\
& $ \ \ \delta  \left ( 3A^2 + k A^2 + t B^2 \right )$ \ if $\varepsilon =-1$,
\\
{\bf VI.}
& \ \ $\delta \left ( (2A^2 -1) + t B^2 \right )$. 
\end{tabular}
\par \medskip
(3) Surgery coefficient: $bB + \delta A$
\par \vskip20pt
\begin{tabular}{l|c|c|c|c}
&
$A$ \qquad $\,$ & $k$ (if $\varepsilon =+1$) 
& $k$ (if $\varepsilon =-1$) 
& $\varepsilon p \not=$
\\
\hline 
{\bf III}
&
\hfill $n+1$ & \hfill $p-1$ & \hfill $p-2$ & $-1,0$
\\
\hline
{\bf IV}
&
\hfill $2n+1$ & \hfill $p-2$ & \hfill $p-3$ & $-2, -1, 0, 1$
\\
\hline
{\bf V}
&
\hfill $2n+3$ & \hfill $p-2$ & \hfill $p-3$ & $-2, -1, 0, 1$
\\
\hline
{\bf VI}
&
\hfill $2n+2$ & - &  \hfill ($0$) & -
\\
\hline
\end{tabular}
\par
\medskip
(4) Prameter Translation (to \cite{Bg, Ba, Ba3, DMM}), 
see Subsection~\ref{sbsec:parameter}.
\par \vskip20pt
\caption{Parameters of Berge's knots}
\label{tbl:Para}
\end{center}
\end{table}
%
%

\pagebreak  

\begin{table}[h]

\begin{tabular}{llll}
{\bf III.} & $\varepsilon =1$, \ $A=2,3,4, \cdots$ \hskip2.3cm $\,$
& & $\varepsilon =-1$, \ $A=2,3,4, \cdots$ \\
& $\textrm{coef} = 6A^2 - 3A $, 
& & $\textrm{coef}  = 6A^2 + 3A$ \\
& $\textrm{area} ({\mathcal L})= 6A^2 -3A$, 
& & $\textrm{area} ({\mathcal L})=6A^2 +3A+1$ \\
\end{tabular}
\begin{center}
\includegraphics[scale=0.5]{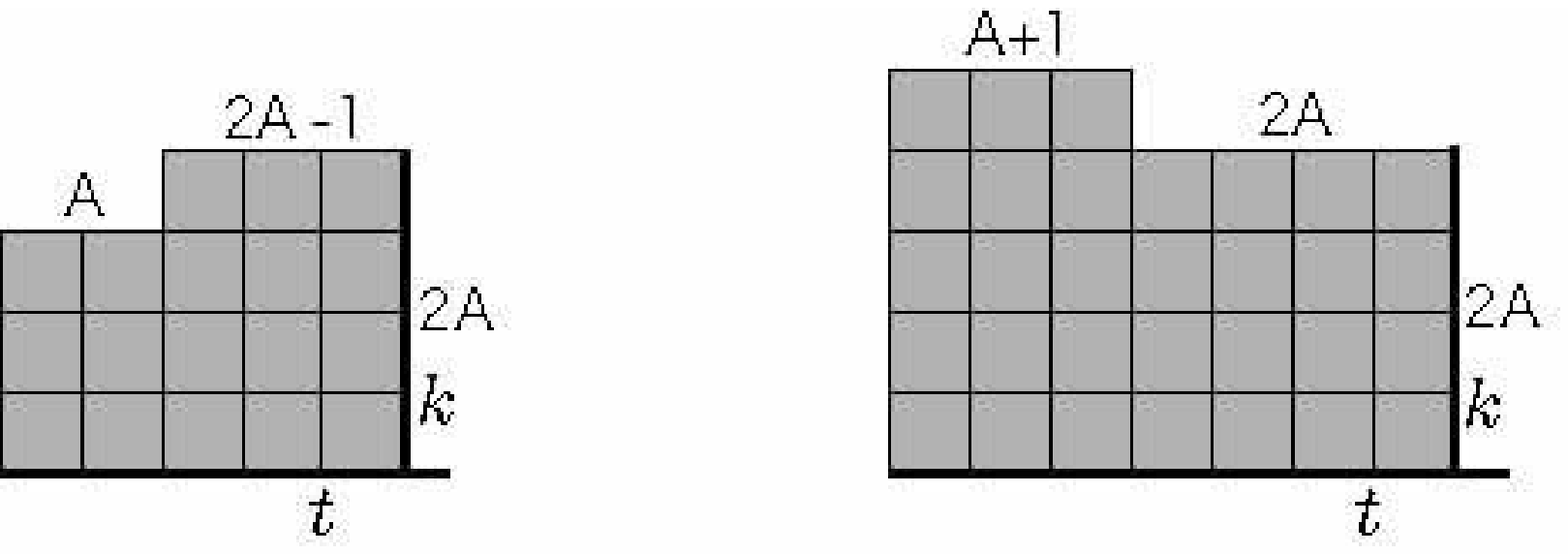}
\end{center}
\begin{tabular}{llll}
{\bf IV.} & $\varepsilon =1$, \ $A=5,7,9, \cdots$\hskip2cm $\,$
& & $\varepsilon =-1$, \ $A=5,7,9, \cdots$\\
& $\textrm{coef} = \frac52 A^2 - \frac32 A$, 
& & $\textrm{coef} = \frac52 A^2 + \frac32 A$ \\
& $\textrm{area} ({\mathcal L}) = \frac52 A^2 - \frac32 A +1 $, 
& & $\textrm{area} ({\mathcal L}) = \frac52 A^2 + \frac32 A $ \\
\end{tabular}
\begin{center}
\includegraphics[scale=0.5]{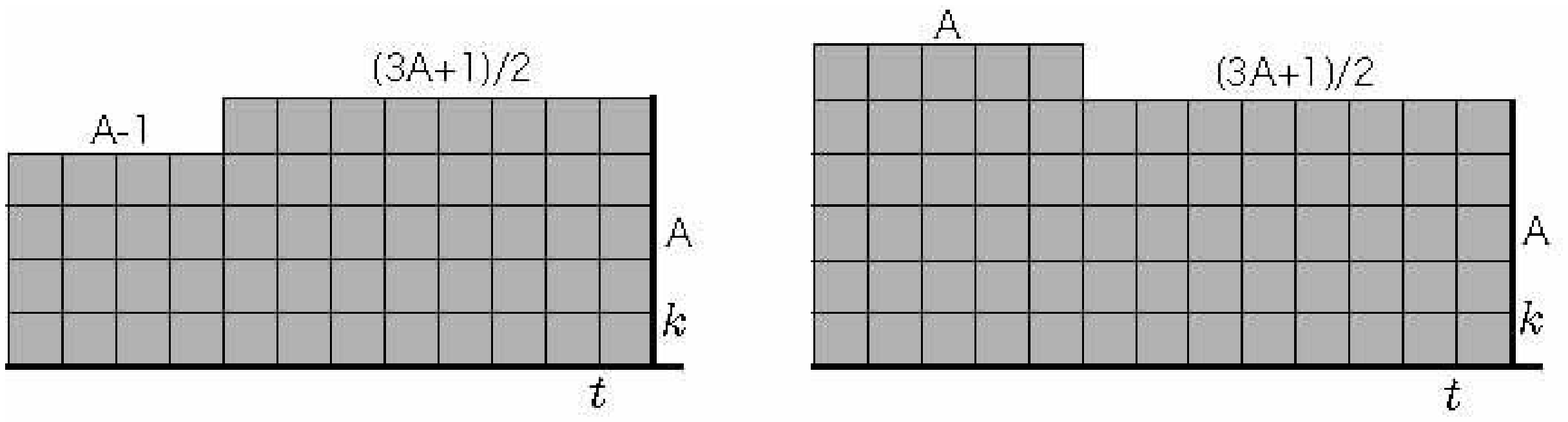}
\end{center}
\begin{tabular}{llll}
{\bf V.} & $\varepsilon =1$, $A= 3,5,7, \cdots$
 \hskip2.3cm $\,$
& & $\varepsilon =-1$, $A= 3,5,7, \cdots$ \\
& $\textrm{coef} = 2A^2$, 
& & $\textrm{coef} = 3A^2$ \\
& $\textrm{area} ({\mathcal L})= 2A^2 +1$,
& & $\textrm{area} ({\mathcal L})= 3A^2$ \\
\end{tabular}
\begin{center}
\includegraphics[scale=0.5]{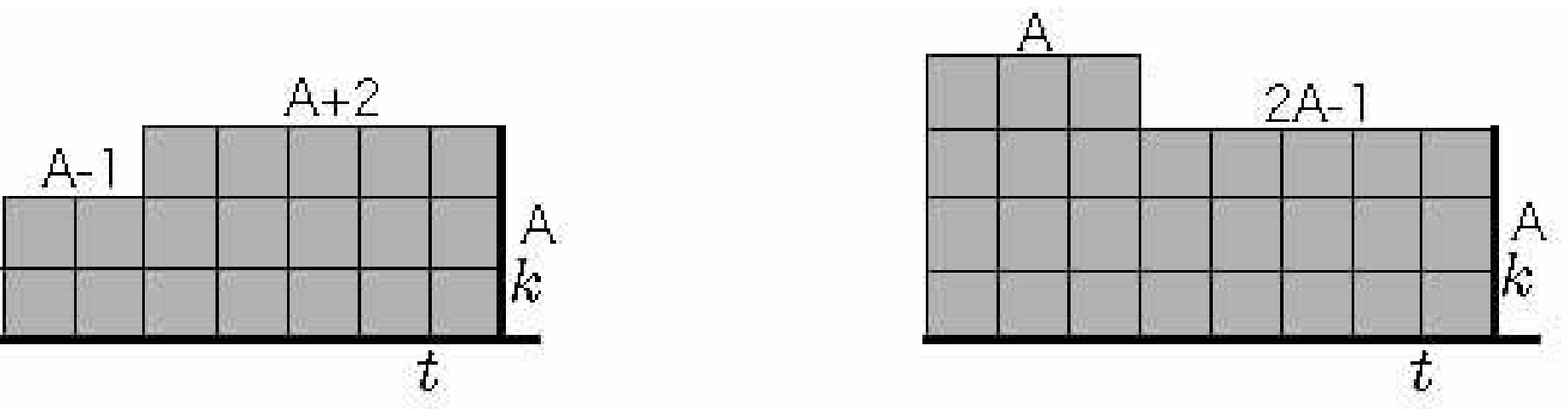}
\end{center}
\begin{tabular}{ll}
{\bf VI.} & $A= 4,6,8, \cdots$ \\
& $\textrm{coef} = 2A^2-1$ \\
& $\textrm{area} ({\mathcal L})= 2A^2$
\end{tabular}\\[-1.5cm]
\begin{center}
\includegraphics[scale=0.5]{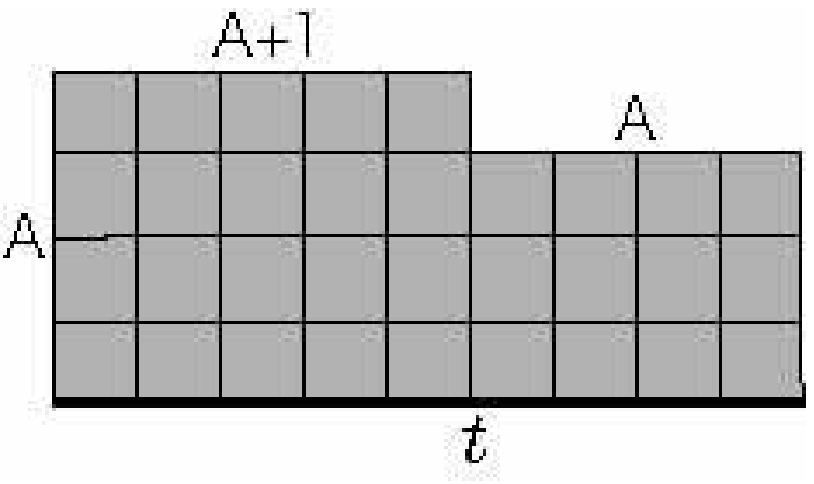}\\
\end{center}
\caption{Berge's knots presented by L-shaped regions}
\label{tbl:Lregion}

\end{table}

\end{document}